\documentclass[11pt]{article}

\usepackage{tikz}
\usepackage{amsmath}
\usepackage{amsthm}
\usepackage{amssymb}
\usepackage{bm}
\usepackage{amscd}
\usepackage{txfonts}
\usepackage{authblk}
\usepackage{amsfonts}
\usepackage{graphicx}
\usepackage{mathtools}
\usepackage{enumitem}
\usepackage{hyperref}

\usepackage[left=3cm,top=3cm,right=3cm,bottom=2cm]{geometry}
\newcommand{\algA}{\mathbb{A}}

\newcommand{\algC}{\mathbb{C}}

\newcommand{\ba}{\mathbf{a}}
\newcommand{\bo}{\mathbf{0}}
\newcommand{\bC}{\mathbf{C}}
\newcommand{\bc}{\mathbf{c}}
\newcommand{\Db}{\Diamondblack} %diamond black
\newcommand{\Bb}{\blacksquare} %box black
\renewcommand{\epsilon}{\varepsilon}
\renewcommand{\i}{\mathbf{i}}
\renewcommand{\j}{\mathbf{j}}

\newcommand{\m}{\mathbf{m}}
\newcommand{\n}{\mathbf{n}}

\renewcommand{\phi}{\varphi}
\newcommand{\mvml}{L_\algA}
\newcommand{\emvml}{L_\algA^+}
\newcommand{\mvfol}{\mathit{FO}_\algA}
\newcommand{\mvfolz}{\mathit{FO}^0_\algA}
\newcommand{\emvfol}{\mathit{FO}_\algA^{+}}

\newcommand{\mvsol}{\mathit{SO}_\algA}
\newcommand{\emvsol}{\mathit{SO}_\algA^{+}}
\newcommand{\Nom}{\mathsf{NOM}}
\newcommand{\Cnom}{\mathsf{CNOM}}
\newcommand{\Prop}{\mathsf{PROP}}
\newcommand{\Var}{\mathsf{VAR}}

\newcommand{\ST}{S\!T}
\newcommand{\modM}{\mathfrak{M}}
\newcommand{\modF}{\mathfrak{F}}
\newcommand{\foM}{\mathcal{M}}
\newcommand{\cfoM}{\mathfrak{M}^{\mathsf{FO}}}

\renewcommand{\P}{\mathbf{P}}

\newcommand{\Pre}{\mathit{PRE}}

\newcommand{\Pos}{\mathsf{POS}}

\newcommand{\Boxat}{\text{\sf BOX-AT}}
\newcommand{\Rel}{\mathsf{REL}}
\newcommand{\oP}{\overline{P}}
\newcommand{\obP}{\overline{\mathbf{P}}}
\renewcommand{\iff}{\text{iff}}

\newcommand{\nomi}{\mathbf{i}}

\newcommand{\cnomm}{\mathbf{m}}

\newcommand{\jty}{J^{\infty}}
\newcommand{\mty}{M^{\infty}}

\newif\ifmargincoms

\margincomstrue
\newtheorem{thm}{Theorem}[section]
\newtheorem{lem}[thm]{Lemma}
\newtheorem{lemma}[thm]{Lemma}
\newtheorem{prop}[thm]{Proposition}

\newtheorem{cor}[thm]{Corollary}

\theoremstyle{definition}
\newtheorem{definition}[thm]{Definition}
\newtheorem{example}[thm]{Example}

\newtheorem{remark}[thm]{Remark}

\theoremstyle{remark}

\title{Correspondence Theory for Many-valued Modal Logic}
\author[1]{Cecelia Britz\thanks{The financial assistance of the 
National Research Foundation (NRF) towards this research is 
hereby acknowledged. Opinions expressed and conclusions arrived 
at, are those of the author and are not necessarily to be 
attributed to the NRF.}}
\author[2]{Willem Conradie}
\author[1]{Wilmari Morton}
\affil[1]{\small Department of Mathematics and Applied Mathematics, 
University of Johannesburg, South Africa}
\affil[2]{School of Mathematics, University of the Witwatersrand, Johannesburg, South Africa}

\begin{document}
\maketitle
\begin{abstract}
The aim of the present paper is to generalise Sahlqvist correspondence theory to the many-valued
modal semantics defined by Fitting, assuming a perfect Heyting algebra as truth value space. We present the standard translations
between many-valued modal languages and suitably defined first-order and second-order
correspondence languages and prove their correctness.  We introduce a notion of many-valued modal frame correspondence with a truth value parameter. Exploring the consequences of this definition, we define many-valued analogues of the syntactically specified classes of Sahlqvist and inductive formulas. We adapt the ALBA algorithm to effectively compute many-valued parameterized local frame correspondents for all many-valued Sahlqvist and inductive formulas. Lastly we prove that the many-valued frame correspondent (parameterized with any non-zero truth value) of every classical Sahlqvist formula is syntactically identical to its standard crisp frame correspondent.  

\end{abstract}

\section{Introduction}

Interest in many-valued modal logics goes back at least to the 1960s, see e.g. \cite{SEGERBERG}, \cite{schotch1978note} and \cite{thomason1978possible}. 
%For example, in \cite{hajek1994modal}, H\'{a}jek and coauthors introduce and study many-valued modal logics of knowledge and belief with $n$-valued Lukasiewicz logic as propositional base. 
In the present paper, we will mostly follow the seminal work of Fitting in \cite{Fitting, Fitting2}, which introduces Kripke semantics for many-valued modal logics in which a finite distributive lattice or Heyting algebra forms the space of possible truth values and formulas are interpreted on generalizations of Kripke models where proportional variables and the accesibility relation take truth values in this algebra.

In the setting of standard (two-valued) modal logic, a modal formula $\phi$ and a first-order formula $\alpha(x)$ in one free variable $x$, are called \emph{local frame correspondents} if $\modF, w \Vdash \phi$ iff $\modF \models \alpha [x:= w]$ for all frames $\modF = (W,R)$ and all states $w \in W$. Most modal formulas lack first-order frame correspondents and the problem of determining which modal formulas have them is undecidable \cite{chagrov1997modal}. Modal correspondence theory has been extensively studied since the 1970s and classical references include \cite{Sahl}, \cite{vanb} and \cite{vanBenthemChapter2001}. 

%\texttt{Refer to correspondence results in \cite{Koutras1}, \cite{Koutras5} and \cite{Frank}.} 
A number of results of correspondence-theoretic interest have appeared in the  literature on many-valued logic. Notable results include those of  Koutras et.\ al.\ in  \cite{Koutras1} and \cite{Koutras5} where properties of many-valued frames defined by the validity of variations on a number of well-known modal axioms are considered. In \cite{Frank} Frankowski considers a filter-parametrized notion of frame validity and the frame properties defined by a number of modal formulas via this kind of validity. More recently, Badia, Caicedo and Noguera \cite{BadiaCaicedoNoguera2023} proved that, under certain assumptions including the finiteness of the truth value algebra, the classes of \emph{crisp} Kripke frames definable by modal formulas under many-valued interpretation coincide with the classes so definable under 2-valued interpretation. 

The aim of the present paper is to generalise Sahlqvist correspondence theory to the many-valued
setting. We settle on a version of Fitting's semantics with an arbitrary perfect (or at least complete) Heyting algebra as the space of truth values.

After collecting the necessary preliminaries in Section \ref{sec:prelims}, we present standard translations
between many-valued modal languages and suitably defined first-order and second-order
correspondence languages in Section \ref{sec:Languages:and:Translations}, and prove their correctness.  In Section \ref{sec:mv:corresp:theory} we formally introduce a notion of many-valued frame correspondence  parametric in a truth value. This definition raises several questions, including: to what extent does the existence and shape of a many-valued first-order correspondent depend on the truth-value parameter? Are there effectively identifiable classes of modal formulas guaranteed to have many-valued first order correspondents (for a given truth value parameter)? If a modal formula has a first-order correspondent, does it have a many-valued first-order correspondent (and vice versa) and, if so, how are these correspondents related? The subsequent sections aim to provide answers to these questions. In Section \ref{sec:Sahl:and:Ind:fmls} we define many-valued analogues of the Sahlqvist and inductive formulas and inequalities (see e.g.\ \cite{vanb,vanBenthemChapter2001}, \cite{Goranko2006-GORECF} and \cite{ALBA, nonDistALBA}) which we illustrate with some examples. In Section \ref{Section:ALBA} we indicate how the ALBA algorithm \cite{ALBA,nonDistALBA} can be adapted to effectively compute truth value-parameterized many-valued frame correspondents of all many-valued Sahlqvist and inductive formulas. This establishes a very general first-order correspondence result for many-valued modal logic. 
%Next we show that algorithmic approaches to obtain first-order local frame correspondents can be generalise to the many-valued setting and present the generalisation of the ALBA algorithm for this purpose.  
Section \ref{Section:SvBResult} is devoted to proving that the parameterized many-valued frame correspondent of every Sahlqvist formula (in the sense of the original classical definition) is syntactically identical to its standard crisp correspondent.  This is done by generalizing the classical Sahlqvist-van Benthem correspondence algorithm to the many-valued setting and proving that its output coincides with that of the classical algorithm. These results falls within a line a research concerned with establishing systematic connection between correspondents for formulas across different relational interpretations: in the present case for modal Sahlqvist formulas interpreted over crisp Kripke frames and many-valued Kripke frames and, e.g. in \cite{conradie2024modal} for the smaller class of Sahlqvist modal reduction principles interpreted over crisp and many-valued Kripke frames, and crisp and many-valued polarity-based frames. We conclude in Section \ref{sec:conclusions}.

\section{Preliminaries}\label{sec:prelims}

Recall that a \emph{Heyting algebra} is an algebra $\algA = (A,\vee,\wedge,\rightarrow,0,1)$ where the reduct $(A,\vee,\wedge,0,1)$ is a bounded distributive lattice and $\rightarrow$ is the right residual of $\wedge$, i.e.\ for all $a, b, c  \in A$ it holds that $a \wedge b \leq c$ iff $a \leq b \rightarrow c$. A Heyting algebra $\algA = (A,\vee,\wedge,\rightarrow,0,1)$ is \emph{complete} if its lattice reduct is, i.e.,\ if $\bigvee S$ and $\bigwedge S$ exist for all $S \subseteq A$. It is easy to check that in any complete Heyting algebra $\bigvee S \rightarrow b = \bigwedge_{s \in S} (s  \rightarrow b)$ and $a \rightarrow \bigwedge S = \bigwedge_{s \in S}(a \rightarrow s)$. Moreover, every complete Heting algebra can be uniquely expanded to a complete \emph{bi-Heyting algebra} $(A,\vee,\wedge,\rightarrow, -, 0,1)$ by setting $a - b = \bigwedge \{c \mid a \leq b \vee c \}$, which makes $-$ the left residual of $\vee$, i.e.,\ for all $a, b, c  \in A$ it holds that $a \leq b \vee c$ iff $a - b \leq c$. The operation $-$ is referred to as the \emph{pseudo difference} or the \emph{co-implication}.

In a complete lattice $\algC$, an element $c\in \algC\setminus \{0\}$
is \emph{completely join-irreducible} iff  $c = \bigvee S$ implies $c \in S$, for every $S\subseteq \algC$. The element $c$ is \emph{completely join-prime}
if, for every $S\subseteq \algC$, if $c\leq \bigvee S$, then $c\leq s$ for some $s\in S$. It is easy to check that all completely join-prime elements are completely join-irreducible and that, moreover, they coincide when $\algC$ is \emph{frame distributive}, i.e.,\ when finite meets distribute over arbitrary joins. The collection of completely join-irreducible elements of $\algC$ will be denoted by $\jty(\algC)$. The completely \emph{meet-irreducible/prime} elements are defined order dually, and the set of completely meet-irreducible elements of $\algC$ is denoted by $\mty(\algC)$.

A complete lattice $\algC$ is {\em perfect} (cf.\ \cite[Def.\ 2.9]{DGP}) if $\jty(\algC)$ is join-dense in $\algC$ and $\mty(\algC)$ is meet-dense in $\algC$, i.e.,\ if for every $u\in \algC$, $u = \bigvee\{r\in \jty(\algC)\ |\ r\leq u\}$ 	and $u = \bigwedge\{z\in \mty(\algC)\ |\ u\leq z\}$. A {\em perfect Heyting algebra} is a Heyting algebra with a perfect lattice reduct. In particular, every finite Heyitng algebra is perfect.

On any perfect distributive lattice $\algC$, we may define two very useful maps, namely  $\kappa: \jty(\algC)\rightarrow \mty(\algC)$
defined by $j \mapsto \bigvee\{u \in \algC \mid j \not\leq u\}$ and its inverse $\lambda: \mty(\algC)\rightarrow \jty(\algC)$, given by the assignment $m \mapsto \bigwedge\{u \in \algC\mid u\not \leq m \}$. These two maps are order isomorphisms between $\jty(\algC)$ and $\mty(\algC)$ considered as subposets of $\algC$ (cf.\ \cite[Sec.\ 2.3]{GEHRKE200565}). It follows from the definitions that, for every $u \in \algC$, every $j \in \jty(\algC)$ and every $m \in \mty(\algC)$,
\begin{equation}\label{eq:kappa}
j \not \leq u\ \mbox{ iff }\ u \leq \kappa(j).
\end{equation}
\begin{equation}\label{eq:lambda}
u \not\leq m \ \mbox{ iff }\ \lambda(m)\leq u.
\end{equation}
For any lattice $\algA$, let $\algA^1=\algA$ and $\algA^\partial$ the lattice with dual ordering.%\wmnote{Do we need (all the detail in) this paragraph?} For an $n$-ary function $f$, its order-type $\varepsilon$ is an $n$-tuple from $\{1,\partial\}^n$. Let $(\mathcal{F},\mathcal{G})$ be a pair of disjoint function symbols with arities given by $\rho:\mathcal{F}\cup\mathcal{G}\to \mathbb{N}$ and order types, such that for $h\in \mathcal{F}\cup\mathcal{G}$, $\varepsilon:h \mapsto\{1,\partial\}^{\rho(h)}$.  A {\em lattice expansion} is a tuple $\mathfrak{L}=(L,\mathcal{F}^\mathfrak{L},\mathcal{G}^\mathfrak{L})$ such that $L$ is a bounded lattice with $\mathcal{F}^\mathfrak{L}=\{f^\mathfrak{L}\mid f\in\mathcal{F}\}$ and $\mathcal{G}^\mathfrak{L}=\{g^\mathfrak{L}\mid g\in\mathcal{G}\}$ such that every $h^\mathfrak{L}\in \mathcal{F}^\mathfrak{L}\cup\mathcal{G}^\mathfrak{L}$ is an $\rho(h)$-ary operation on $L$ of order-type $\varepsilon(h)$. If $L$ is distributive (resp. perfect), $\mathfrak{L}$ is called a distributive (resp. perfect) lattice expansion.  Moreover, $\mathfrak{L}$ is called {\em normal} if every $f^\mathfrak{L}\in \mathcal{F}^\mathfrak{L}$ (resp. $g^\mathfrak{L}\in \mathcal{G}^\mathfrak{L}$), if  $\varepsilon(f)$'s $i^{\text{th}}$ coordinate is $1$  (resp. $\varepsilon(g)$'s $i^{\text{th}}$ coordinate is $1$), then $f^\mathfrak{L}$ preserves finite (potentially empty) joins in the $i^{\text{th}}$ coordinate  (resp. $g^\mathfrak{L}$ preserves finite (potentially empty) meets in the $i^{\text{th}}$ coordinate); and if $\varepsilon(f)$'s $i^{\text{th}}$ coordinate is $\partial$ (resp. $\varepsilon(g)$'s $i^{\text{th}}$ coordinate is $\partial$), then $f^\mathfrak{L}$ reverses finite (potentially empty) meets in the $i^{\text{th}}$ coordinate  (resp. $g^\mathfrak{L}$ reverses finite (potentially empty) joins in the $i^{\text{th}}$ coordinate).

%%%%%%
Let $\algA=(A,\vee,\wedge,\rightarrow,0,1)$ be a perfect
Heyting algebra and $\Prop$ a countably infinite set of proposition
letters.  The formulas of the basic many-valued modal language over
$\algA$ (see \cite{Fitting}), $\mvml$, are given by the following recursive definition:
\begin{align*}
	\phi := \ba \mid p \mid \phi \vee \psi \mid \phi \wedge \psi 
	\mid \phi \rightarrow \psi \mid \Box \phi \mid \Diamond \phi
\end{align*}
where $p\in\Prop$ and each $\ba$ is a propositional constant that corresponds to an $a\in A$.
Negation can be defined in terms of $\rightarrow$ and $\bo$ as $\neg \phi := \phi \rightarrow \bo$. An expression of the form $\phi \leq \psi$, where $\phi$ and $\psi$ are $\mvml$-formulas, is called an {\em $\mvml$-inequality}. An expression of the form $\phi_1 \leq \psi_1 \& \cdots \& \phi_n \leq \psi_n \Rightarrow \phi_0 \leq \psi_0$ where the $\phi_i$ and $\psi_i$, $0 \leq i \leq n$, are $\mvml$-formulas is called an {\em $\mvml$-quasi-inequality}.

A {\em(Kripke) frame for $\mvml$} ({\em $\algA$-frame}) is a pair $\modF = (W,R)$ such that $W$ is nonempty and $R$ is a binary $\algA$-fuzzy relation
%The term fuzzy relation is used in Hajek's book, eg. Chapter five page 109 line -2
on $W$, i.e., $R:(W\times W)\rightarrow A$. A {\em valuation on $\modF$} is a binary function $V: \left(\Prop\times W\right)\rightarrow A$ that assigns to every pair $(p,w)$ the truth-value of $p$ at $w$.
A {\em (Kripke) model for $\mvml$} ({\em $\algA$-model}) is a pair $\modM=(\modF, V)$
such that $\modF$ is an $\algA$-frame and $V$ a valuation on $\modF$.
A valuation $V$ can be extended to all formulas of $\mvml$ (also denoted by $V$) in the following way:
\begin{enumerate}[label=(\roman*)]
	\item $V(\ba,w)=a$ for all $a\in A$,
	\item $V(\phi \star \psi,w) = V(\phi,w) \star^{\algA} V(\psi,w)$ for $\star\in\{\vee,\wedge,\rightarrow\}$,
	\item $V(\Diamond\phi, w) = \bigvee\left\{R(w,u) \wedge^{\algA} V(\phi,u)\mid u\in W\right\}$, and
	\item $V(\Box\phi, w) = \bigwedge \left\{R(w,u) \rightarrow^{\algA} V(\phi,u) \mid u\in W\right\}$.
\end{enumerate}
Note that, since $\algA$ is complete, $V(\Diamond \psi, w)$ and $V(\Box \psi, w)$ will always exist in $\algA$.

Notions of truth and validity in $\algA$-frames $\modF=(W,R)$ and $\algA$-models $\modM=(\modF,V)$, parametrized with truth values from $\algA$, are defined following \cite{Koutras}. That is, if $a\in A$, $w\in W$ and $\phi$ a formula of $\mvml$, then $\phi$ is said to be {\em $a$-true at $w$ in $\modM$} (notation: $\modM,w\Vdash_{a} \phi$)
if $V(\phi,w)\geq a$, and  \emph{true at $w \in \modM$} (notation: $\modM,w\Vdash \phi$) if $\modM,w\Vdash_{1} \phi$.  A formula $\phi$ of $L_{\algA}$ is 
{\em $a$-true in $\modM$} (notation: $\modM\Vdash_{a} \phi$) if $\modM,w\Vdash_{a} \phi$ for all $w\in W$ and \emph{true $\modM$} (notation: $\modM\Vdash \phi$) if $\modM\Vdash_{1} \phi$.

An $\mvml$-inequality $\phi \leq \psi$ is {\em true in $\modM$ at $w$} (notation: $\modM,w\Vdash \phi \leq \psi$) if $V(\phi, w) \leq^{\algA} V(\psi, w)$, and  $\phi \leq \psi$ is {\em true in $\modM$} (notation: $\modM \Vdash \phi \leq \psi$) if $V(\phi, w) \leq^{\algA} V(\psi, w)$ for all $w \in W$. 

We note the following further reductions between the notions above:  
\begin{equation}
\begin{aligned}
    \modM, w \Vdash \phi \leq \psi &\quad \text{iff} \quad \modM, w \Vdash \phi \to \psi &\modM \Vdash \phi \leq \psi &\quad \text{iff} \quad \modM \Vdash \phi \to \psi\\ 
    \modM, w \Vdash_{a} \phi \to \psi &\quad \text{iff} \quad \modM, w \Vdash a \wedge \phi \to \psi &\modM \Vdash_{a} \phi \to \psi &\quad \text{iff} \quad \modM \Vdash a \wedge \phi \to \psi\\
    \modM, w \Vdash_{a} \phi  &\quad \text{iff} \quad \modM, w \Vdash a \leq \psi &\modM \Vdash_{a} \phi  &\quad \text{iff} \quad \modM \Vdash a \leq \psi\\ 
    \modM, w \Vdash a \leq \psi &\quad \text{iff} \quad \modM, w \Vdash a \to \psi &\modM \Vdash a \leq \psi &\quad \text{iff} \quad \modM \Vdash a \to \psi
\end{aligned}
\label{eq:satisfac:eqivs}
\end{equation}
This suggests that we can define the \emph{$a$-truth of an inequality} by declaring that $\modM,w \Vdash_{a} \phi \leq \psi$ if $\modM, w \Vdash a \wedge \phi \leq \psi$, and that $\modM \Vdash_{a} \phi \leq \psi$ if $\modM \Vdash a \wedge \phi \leq \psi$. 

An $\mvml$-quasi-inequality $\phi_1 \leq \psi_1 \& \cdots \& \phi_n \leq \psi_n \Rightarrow \phi_0 \leq \psi_0$ is {\em true in $\modM$} (notation: $\modM \Vdash \phi_1 \leq \psi_1 \& \cdots \& \phi_n \leq \psi_n \Rightarrow \phi_0 \leq \psi_0$) if $\modM \Vdash \phi_0 \leq \psi_0$ or $\modM \not \Vdash \phi_i \leq \psi_i$ for some $1 \leq i \leq n$.

If $a\in A$, $w\in W$ and $\phi$ a 
formula of $\mvml$, then $\phi$ is said to be {\em $a$-valid in $\modF$ at $w$} 
(notation: $\modF,w \Vdash_{a} \phi$) if $V(\phi,w)\geq a$ for all valuations $V$ on $\modF$.
A formula $\phi$ of $L_{\algA}$ is {\em $a$-valid} in an $\algA$-frame $\modF$ 
(notation: $\modF\Vdash_{a} \phi$) if $\modF,w\Vdash_{a} \phi$ for all $w\in W$. 
An $\mvml$-inequality $\phi \leq \psi$ is {\em valid in $\modF$ at $w$} (notation: $\modF,w\Vdash \phi \leq \psi$) if $V(\phi, w) \leq^{\algA} V(\psi, w)$ for all valuations $V$ on $\modF$, and  $\phi \leq \psi$ is {\em valid in $\modF$} (notation: $\modF \Vdash \phi \leq \psi$) if $\modF,w\Vdash \phi \leq \psi$ for all $w \in W$. An $\mvml$-quasi-inequality $\phi_1 \leq \psi_1 \& \cdots \& \phi_n \leq \psi_n \rightarrow \phi_0 \leq \psi_0$ is {\em valid in $\modF$} (notation: $\modF \Vdash \phi_1 \leq \psi_1 \& \cdots \& \phi_n \leq \psi_n \rightarrow \phi_0 \leq \psi_0$) if it is true in all models based on $\modF$.

The {\em complex algebra} $\modF^+$ of an $\algA$-frame $\modF=(W,R)$  is the algebra
%	\[
$\modF^+ = (\algA^W,\Diamond_{R}, \Box_{R},\{\ba\}_{a\in A})$,
%	\]
where $\algA^W=\{f\mid f:W\to \algA\}$ (i.e., the set of functions from $W$ into $\algA$),
$\ba: W \rightarrow A$ is the constant function $\ba(w)=a$ for all $a\in \algA$ and $w\in W$, and
\begin{align*}
(\Diamond_{R} f)(w) &= \bigvee \left\{f(u) \wedge R(w,u)\mid u\in W\right\}\\
(\Box_{R} f)(w) &= \bigwedge \left\{R(w,u)\rightarrow f(u)\mid u\in W \right\}.
\end{align*} 

\begin{lem}\label{lem:complex:is:perf:Heyting}
	Let $\mathfrak{F}$ be an $\algA$-frame. Then $\mathfrak{F}^+ = (\algA^W,\Diamond_{R}, \Box_{R},\{\ba\}_{a\in A})$ is a perfect Heyting-algebra with additional operators $\Diamond_{R}$ and $\Box_{R}$ that are, respectively, completely join-preserving and completely meet-preserving.
\end{lem}
\begin{proof}
	Since the class of Heyting algebras is a variety, it is closed
	closed under products, so $\algA^W$ is a Heyting algebra.  Its completeness and complete distributivity
	follows from $\algA$'s completeness and complete distributivity and the fact that $\vee$ and $\wedge$ are defined coordinate-wise
	on $\algA^W$.	
	
	Let $\{f_i \mid i \in I\} \subseteq \algA^W$. Since $\vee$ in $\mathfrak{F}^+$ is defined coordinate-wise, it is enough to show that   $(\Diamond_{R}\left( \bigvee_{i \in I} f_i \right))(w) = \bigvee_{i \in I} [(\Diamond_{R}\left(  f_i \right))(w)]$ for all $w \in W$. Indeed, $(\Diamond_{R}( \bigvee_{i \in I} f_i))(w) = \bigvee_{u\in W} \left( (\bigvee_{i \in I} f_i)(u) \wedge R(w,u)\right) = \bigvee_{u\in W} \left(\bigvee_{i \in I}f_i(u) \wedge R(w,u) \right) = \bigvee_{i \in I} \bigvee_{u\in W} \left(f_i(u) \wedge R(w,u) \right) = \bigvee_{i \in I} [(\Diamond_{R}(f_i))(w)]$.
	
	The case for $\Box_{R}$ is similar. 
\end{proof}
\noindent Formulas of $\mvml$ can be interpreted as algebraic terms in Heyting algebras with unary operators $\Diamond$ and $\Box$ in the standard way as, hence, can  $\mvml$-(in)equalities and quasi-(in)equalities. In the case of $\algA$-frames and their complex algebras, the relational and algebraic semantics are connected in the expected way, as is made explicit in the following, routine lemma:
\begin{lem}\label{lem:valid:in:frame:and:complex:alg}
For any $\algA$-frame $\modF=(W,R)$, $\mvml$-inequality $\phi \leq \psi$ and $\mvml$-quasi-inequality $\phi_1 \leq \psi_1 \& \cdots \& \phi_n \leq \psi_n \rightarrow \phi_0 \leq \psi_0$, 
\[
\modF \Vdash \phi \leq \psi \quad \text{ iff } \quad \modF^+ \models \phi \leq \psi
\]
and
\[
\modF \Vdash \phi_1 \leq \psi_1 \& \cdots \& \phi_n \leq \psi_n \rightarrow \phi_0 \leq \psi_0 \quad \text{ iff } \quad \modF^+ \models \phi_1 \leq \psi_1 \& \cdots \& \phi_n \leq \psi_n \rightarrow \phi_0 \leq \psi_0.
\]
\end{lem}

\section{Many-valued languages and translating between them}\label{sec:Languages:and:Translations}

In this section we introduce several languages that we will need for our correspondence-theoretic investigations. We begin in subsection \ref{sec:ext:ml:lang} by formulating an extension $\emvml$ of $\mvml$ which we will use to capture correspondence-theoretic arguments and is the language of the ALBA calculus (see Section \ref{Section:ALBA}). Subsection \ref{sec:coresp;languages} introduces the first and second-order model and frame correspondence languages into which $\mvml$ and $\emvml$ can be faithfully translated by the standard translations introduced and proven correct in Subsection \ref{Sec:StndTrsltn}.

%\wmnote{Refs to literature. Why Heyting and not residuated lattice?}

\subsection{The extended modal language}\label{sec:ext:ml:lang}

%The extended many-valued modal language over A and its Kripke semantics
In order to perform correspondence theoretic reductions in the ALBA calculus (Section \ref{Section:ALBA}) we will need to enhance the expressive power of $\mvml$ in three directions: by adding inverse (or tense) modalities, a connective $-$ to be interpreted with the co-implication of $\algA$, as well as special variables called nominals and co-nominals, analogous to nominals in hybrid logic. Formally, given disjoint and countably infinite sets $\Prop$, $\Nom$ and $\Cnom$ of propositional variables, nominals and co-nominals, respectively, the formulas of the language $\emvml$ are given by   
\[
\phi := \ba \mid \i \mid \m \mid p \mid \phi \vee \psi \mid \phi \wedge \psi 
\mid \phi \rightarrow \psi\mid \phi - \psi \mid \Box \phi \mid \Diamond \phi \mid \Db \phi \mid \Bb \phi
\label{indformplus}
\]
where $\i \in \Nom$, $\m \in \Cnom$, $p \in \Prop$  and each $\ba$ is a propositional constant that corresponds to an $a\in A$.
  
$\emvml$ is interpreted over $\algA$-frames and $\algA$-models, but where a valuation $V$ on an $\algA$-model is now a binary function $V:((\Prop\cup\Nom\cup\Cnom) \times W)\rightarrow A$ that additionally satisfies: 
\begin{enumerate}[label=(\roman*)]
\item $V(\i,w)\in J^{\infty}(\algA)$ for exactly one $w\in W$ and $V(\i,u)=0$ for all $u\in W$ such that $w \neq u$.
\item $V(\m,w)\in M^{\infty}(\algA)$ for exactly one $w\in W$ and $V(\m,u)=1$ for all $u\in W$ such that $w \neq u$.
\end{enumerate}
Thus a nominal (co-nominal) takes a non-$0$ (non-$1$) value at exactly one state in a model and this value must be a join (meet) irreducible element of $\algA$.  This interpretation of nominals and co-nominals is motivated by the fact that we wish them to range over join- and meet-irreducibles of the complex algebra of $\algA$-frames, respectively. The additional clauses needed to extend the domain of a valuation $V$ to all $\emvml$-formulas are:
\begin{enumerate}[label=(\roman*)]
\item $V(\phi - \psi,w) = V(\phi,w) -^{\algA} V(\psi,w)$.
\item $V(\Db\phi, w) = \bigvee\left\{ R(u,w) \wedge^{\algA} V(\phi,u)\mid u\in W\right\}$.
\item $V(\Bb\phi, w) = \bigwedge\left\{ R(u,w) \rightarrow^{\algA} V(\phi,u) \mid u\in W\right\}$.
\end{enumerate}
Notions of local and global $a$-truth and validity are extended to $\emvml$ in the obvious way.

%We will also have use for \emph{inequalities} and \emph{quasi-inequalities}.
%\begin{definition}
%\normalfont{
%Given formulas $\phi,\phi_{1},\hdots,\phi_{n}$ and $\psi,\psi_{1},\hdots,\psi_{n}$ from $L^{+}$, a \emph{quasi-inequality} is an expression of the form $(\phi_{1}\leq\psi_{1}\&\hdots\&\phi_{n}\leq\psi_{n})\Rightarrow (\phi\leq\psi)$.  Each expression $\phi_{i}\leq\psi_{i}$ will be referred to as an \emph{inequality}.  
%}
%\label{QuasiI}
%\end{definition}  
%$V(\phi\leq\psi,w)=V(\phi,w)\leq V(\psi,w)$.
%$V\Big((\phi_{1}\leq\psi_{1} \& \hdots \& \phi_{n}\leq\psi_{n})\Rightarrow(\phi\leq\psi),w\Big)=\Big(V(\phi_{1}\leq\psi_{1},w) \wedge^{\algA} \hdots \wedge^{\algA} V(\phi_{n}\leq\psi_{n},w)\Big)\rightarrow^{\algA} V(\phi\leq\psi,w)$.
%Strictly speaking, a valuation applied to a quasi-inequality should interpret $\Rightarrow$ classically, but since an inequality can only be $0$ or $1$ in the current setting, $\rightarrow^{\algA}$ will only receive classical values and hence will behave classically.
%\begin{remark}
%\normalfont{With each nominal $\i$ we can associate a function $f_{\i}\in A^{W}$ such that $f_{\i}(w_{0})\in J^{\infty}(\algA)$ and $f_{\i}(w)=0$ for all $w\in W$ such that $w_{0} \neq w$.  In other words, $V(\i,w_{0})=f_{\i}(w_{0})$.  A similar function $f_{\m}\in A^{W}$ can be defined for each conominal $\m$.  
%}
%\label{NominalF}
%\end{remark}

\subsection{Correspondence languages}\label{sec:coresp;languages}
%Seen through a predicate-logic lens,  
$\algA$-models and frames naturally support the interpretation of first-and second-order languages. These are the correspondence languages for $\mvml$. Because of their importance for the purposes of this paper, we take some time to introduce them  and their semantics formally. For an extended treatment of many-valued predicate logics, we refer the reader to \cite[Chapter 5]{hajek1998metamathematics}.

\subsubsection{$\mvfol$, the first-order correspondence language of $\mvml$, and its semantics}
The \emph{first-order correspondence language of $\mvml$}, denoted $\mvfol$, is the first-order language with equality built over a denumerable set of individual variables $\Var$, a constant symbol $\ba$ for each $a\in A$, a unary predicate symbol $\P$ for each $p\in \Prop$ and a binary relation symbol $R$, using connectives $\vee,\wedge,\rightarrow$ and quantifiers $\forall$ and $\exists$. As usual, the unary connective $\neg$ will be an abbreviation for $\neg \alpha:=\alpha\rightarrow \mathbf{0}$ where $\alpha$ is some formula in $\mvfol$. The set of terms of $\mvfol$ coincides  with $\Var$. The \emph{first-order frame correspondence language of $\mvml$}, denoted $\mvfolz$, is obtained from $\mvfol$ by omitting the predicate symbols $P$.

An \emph{$\algA$-model for $\mvfol$} is a tuple $\foM=(W,R,\{P\}_{p\in\Prop})$ with 
a non-empty domain $W$, a binary $\algA$-fuzzy relation $R :W\times W\rightarrow A$, and an $\algA$-fuzzy set $P:W\rightarrow A$ for each unary predicate symbol $\P$. An {\em assignment} on $\foM$ is a function $v: \Var \rightarrow W$
that assigns an element from $W$ to each individual variable $x\in \Var$.
Two assignments $v$ and $v'$ on $\foM$ are said to be {\em  $x$-variants}, denoted $v'\sim_{x} v$, if $v(y)=v'(y)$ for all $y\neq x$.
The truth-value $\left\Vert \alpha\right\Vert_{\foM,v}$ of a formula $\alpha$
of $\mvfol$ in $\foM$ under an assignment $v$ is defined as follows:
for $t,t_1,t_2$ terms of $\mvfol$ and $\alpha,\beta$ formulas of $\mvfol$
\begin{enumerate}[label=(\roman*)]
\item $\Vert t_{1}= t_{2}\Vert_{\foM,v} = 	
\begin{cases}
		1 & \text{if } v(t_{1})=v(t_{2}), \\
		0 & \text{if } v(t_{1})\neq v(t_{2}).\\
\end{cases}$
\item $\Vert R(t_1,t_2)\Vert_{\foM,v}=R\left(v(t_1),v(t_2)\right)$.
\item $\Vert \P(t)\Vert_{\foM,v}=P\left(v(t)\right)$.
\item $\Vert\ba\Vert_{\foM,v}=a$.
\item $\Vert \alpha\star \beta\Vert_{\foM,v}=\Vert\alpha\Vert_{\foM,v}
\star^{\algA} \Vert\beta\Vert_{\foM,v}$ where $\star\in\{\vee,\wedge,\rightarrow\}$.
\item $\Vert\forall x\,\alpha\Vert_{\foM,v}
=\bigwedge\left\{\Vert\alpha\Vert_{\foM,v'}\mid v'\sim_{x} v\right\}$.
\item $\Vert\exists x\,\alpha\Vert_{\foM,v}
=\bigvee\left\{\Vert\alpha\Vert_{\foM,v'} \mid v'\sim_{x} v\right\}$.
\end{enumerate}

As in the classical case, the truth value of a sentence $\alpha$ is independent of any particular assignment, so writing $\Vert\alpha\Vert_{\foM}$ is unambiguous, if $\alpha$ is a sentence. An \emph{$\algA$-model for $\mvfolz$} is any reduct of an $\algA$-model for $\mvfol$ obtained by omitting the fuzzy predicates $P$, i.e., $\algA$-models for $\mvfolz$ are identical to $\algA$-frames for $\mvml$.

Let $a\in A$.  A formula $\alpha$ of $\mvfol$ is {\em $a$-true} 
in $\foM$, denoted by $\foM \models_{a} \alpha$, if $\Vert\alpha\Vert_{\foM,v}\geq a$
for all assignments $v$ on $\foM$. 
If $x\in \Var$ and $w \in W$, then $\alpha$ is
{\em $a$-true under substitution of $w$ for $x$ in $\alpha$} in $\foM$,
denoted by $\foM \models_{a} \alpha[x \coloneqq w]$, 
if $\Vert\alpha\Vert_{\foM,v}\geq a$ for all assignments $v$ on $\foM$ such that $v(x)=w$.
If $\algA$, $\foM$ and $v$ are clear from the context, we will write $\Vert\alpha\Vert$ instead of $\Vert\alpha\Vert_{\foM,v}$. The notion of $a$-truth (under substitution) of $\mvfolz$-formulas in $\algA$-models for $\mvfolz$ is defined in the obvious, analogous way. 

As we will see in Section \ref{Sec:StndTrsltn}, $\mvml$ can be faithfully translated into $\mvfol$. In order to  extend this translation to $\emvml$, we extend $\mvfol$ to $\emvfol$: 
\subsubsection{$\emvfol$, the first-order correspondence language of $\emvml$, and its semantics}

The language $\emvfol$ extend $\mvfol$ by adding:

\begin{enumerate}[label=(\roman*)]
\item an individual constant $\bc_\i$ for each $\i\in\Nom$ (intended to name the state picked out by $\i$) and an individual constant $\bc_\m$for each $\m\in\Cnom$ (intended to name the state singled out by $\m$),
\item a new {\em truth-value} constant symbol $\bC_\i$ for each $\i\in\Nom$ (to record the truth value of $\i$ at $\bc_\i$) and 
a new {\em truth-value} constant symbol $\bC_\m$ for each $\m\in\Cnom$ (to record the truth value of $\m$ at $\bc_\m$), and
\item the binary connectives $-$ and $\preceq$.
%
%\marginpar{Cecelia included "a binary relation symbol $\preceq$ and", but I'm not sure this should be part of the language.} W: Agree. W (2024-07-18): No, now I see that we need it for the traslation of quasi-inequalities.  
\end{enumerate}
Note that the set terms of $\emvfol$ extends that of $\mvfol$ with the constants $\bc_\i$ and $\bc_\m$ and that all $\bC_\i$ and $\bC_\m$ are atomic formulas of $\emvfol$.

A tuple
$\foM=(W,R,\{P\}_{p\in\Prop},\{c_\i\}_{\i\in\Nom},\{C_\i\}_{\i\in\Nom},
\{c_\m\}_{\m\in\Cnom},\{C_\m\}_{\m\in\Cnom})$
is an $\algA$-model for $\emvfol$ if 
\begin{enumerate}[label=(\roman*)]
\item $\foM=(W,R,\{P\}_{p\in\Prop})$ is an $\algA$-model for $\mvfol$,
\item $c_\i\in W,C_\i\in \jty(\algA)$ for each $\i\in\Nom$, and
\item $c_\m\in W,C_\m\in \mty(\algA)$ for each $\m\in\Cnom$.
\end{enumerate}
The interpretation of terms is extended to $\emvfol$ by stipulating that $\Vert \bc_i \Vert_{\foM,v} = c_\i$ and $\Vert \bc_m \Vert_{\foM,v} = c_\m$. The definition of the value of terms and the truth value of formulas in a model $\foM$
under an assignment $v$ is extended to $\emvfol$ by adding the following clauses to those for $\mvfol$: 
\begin{enumerate}[label=(\roman*)]
\item $\Vert \bc_\i \Vert_{\foM,v} = c_{\i}$
\item $\Vert \bc_\m \Vert_{\foM,v} = c_{\m}$
\item $\Vert\bC_{\i}\Vert_{\foM,v}=C_{\i}$ and $\Vert\bC_{\m}\Vert_{\foM,v}=C_{\m}$
\item $\Vert\alpha-\beta\Vert_{\foM,v}
=\Vert\alpha\Vert_{\foM,v}-^{\algA}\Vert\beta\Vert_{\foM,v}$
\item $\Vert\alpha \preceq \beta \Vert_{\foM,v}
= \begin{cases}
		1 & \text{if } \Vert\alpha\Vert_{\foM,v} \leq^{\algA} \Vert\beta\Vert_{\foM,v}, \\
		0 & \text{otherwise}\\
\end{cases}$ 
\end{enumerate}
The notion of \emph{$a$-truth (under substitution) in $\algA$-models} is extended to $\emvfol$  in the expected way. 

\subsubsection{$\mvsol$, the second-order correspondence language of $\mvml$, and its semantics}

The \emph{second-order correspondence language of $\mvml$}, denoted $\mvsol$, has the same non-logical symbols as $\mvfol$ except that the unary predicate symbols $P$ associated with the proposition letters from $\mvml$ are now regarded as \emph{predicate variables}. The set of formulas of $\mvsol$ extends that of $\mvfol$ by closing under the following formation rule: if $\alpha$ is a formula of $\mvsol$ and $P$ is a predicate variable, then  $\forall P \alpha$ and $\exists P \alpha$ are formulas of $\mvsol$. Note that $R$ remains a predicate symbol, i.e. it does not become a variable.  

An \emph{$\algA$-model for $\mvsol$} is just an $\algA$-frame $\modF=(W,R)$.  An {\em assignment for $\mvsol$} on $\modF$ assigns an element from $W$ to each individual variable $x\in \Var$ and an $\algA$-fuzzy set to each predicate variable $P$. The notion of a $P$-variant of an assignment (notation: $v'\sim_P v$) is defined in the obvious way. The truth-value $\left\Vert \alpha\right\Vert_{\modF,v}$ of a formula $\alpha$ of $\mvsol$ in $\modF$ under an assignment $v$ is defined as follows: 
\begin{enumerate}[label=(\roman*)]
%\item If $\alpha$ has no occurrences of predicate variables, then $\alpha$ is a first-order formula and $\Vert\alpha\Vert_{\modF,v}$ is calculated as for $\mvfol$ above
\item all clauses other than those for atomic formulas $P(t)$ and predicate quantifiers are the same as for $\mvfol$. 
\item $\Vert P(t)\Vert_{\modF,v} = v(P)\left(v(t)\right)$.
\item $\Vert\forall P\alpha\Vert_{\modF,v}=
\bigwedge \left\{\Vert \alpha\Vert_{\modF,v} \mid  v'\sim_P v\right\}$.
\item $\Vert\exists P\alpha\Vert_{\modF,v}=
\bigvee \left\{\Vert \alpha\Vert_{\modF,v} \mid  v'\sim_P v\right\}$. 
\end{enumerate}

\subsubsection{$\emvsol$, the second-order correspondence language of $\emvml$, and its semantics}

The extended many-valued second-order language over $\algA$, $\emvsol$, is the same as $\emvfol$, except that the parts of language used to translate nominals and co-nominals now need to be variables rather than constants. 
%but the additional individual constants in $\nCon$ and $\cCon$ and truth value constants are now treated as variables --- specifically, in this context, the elements of $\nCon$ and $\cCon$ are to be regarded as individual constants and the truth value constants $\bC_{\i}$ and $\bC_{\m}$ as nullary predicate variables. 

More precisely, in $\emvsol$ the symbols $\bc_\i$ for each $\i\in\Nom$  and $\bc_\m$ for each $\m\in\Cnom$ (which were individual constants in $\emvfol$) are regarded as individual variables, and the symbols  $\bC_\i$ for each $\i\in\Nom$ and $\bC_\m$ for each $\m\in\Cnom$ (which were truth value constants in $\emvfol$) are as nullary predicate variables.

%Formally, the sets of additional individual constants in $\emvfol$, $\nCon$ and $\cCon$, are replaced by two disjoint sets of individual variables, $\nVar$ and $\cVar$ (both also disjoint form $\Var$)  such that for each $\i\in\Nom$ there is an individual variable $\bc_\i\in\nVar$ and for each $\m\in\Cnom$ there is an individual variable $\bc_\m\in\cVar$. In addition, $\nCon$ and $\cCon$, are replaced by two disjoint sets of individual variables, $\nVar$ and $\cVar$ (both also disjoint form $\Var$)  such that for each $\i\in\Nom$ there is an individual variable $\bc_\i\in\nVar$ and for each $\m\in\Cnom$ there is an individual variable $\bc_\m\in\cVar$.
%It also includes two disjoint sets of nullary predicate variables, $\nPred$ and $\cPred$  (both also disjoint from $\Pred$) to replace the truth-value constants contained in $\emvfol$.  
%
%Note that quantifying over an element $C_\i\in\nPred$ or $C_\m\in\cPred$ provides a way in $\emvsol$ to quantify over all truth-values, that is, the elements of $\algA$.  
%

An {\em assignment for $\emvsol$} on an $\algA$-frame $\modF=(W,R)$ is an assignment for $\mvsol$ extended to assign an element of $W$ to each individual variable $\bc_\i$ and $\bc_\m$ as well as an element 
of $\jty(\algA)$ to each nullary predicate variable $\bC_\i$ and an element of $\mty(\algA)$ to each nullary predicate variable $\bC_\m$. For any nullary predicate variable $\bC$ (either a $\bC_\i$ or a $\bC_\m$ ) a \emph{$\bC$-variant} of an assignment $v$ (notation: $v'\sim_{\bC} v$) is defined in the obvious way. Given a nullary predicate variable $\bC$, $\alpha$ an  $\emvsol$-formula, $\modF=(W,R)$ an $\algA$-model and $v$ an assignment on $\modF$, the clauses for evaluation of quantification over $\bC$ are as expected: 
\begin{enumerate}[label=(\roman*)]
\item $\Vert\forall \bC\alpha\Vert_{\modF,v}=
\bigwedge \left\{\Vert\alpha\Vert_{\modF,v'} \mid 
v'\sim_{\bC} v\right\}$.
\item $\Vert\exists \bC\alpha\Vert_{\modF,v}=
\bigvee \left\{\Vert\alpha\Vert_{\modF,v'} \mid 
v'\sim_{\bC} v\right\}$.
\end{enumerate}

The notion of \emph{$a$-truth (under substitution) in $\algA$-models} is extended to $\mvsol$ and $\emvsol$ in the expected way. 

\subsection{The standard translation}\label{Sec:StndTrsltn}

It is well known that standard modal languages correspond to fragments first- and second-order logic when interpreted on models and frames, respectively. In this section we show that the same is true for many-valued modal logic with respect to the first-order and second-order correspondence languages introduced above. We present the standard translation and show that it is correct.

\begin{definition}
Let $x$ be a first-order individual variable.  
The standard translation $\ST_{x}$ of $\emvml$ into $\emvfol$ is given by the following clauses:
\begin{enumerate}[label=(\roman*)]
\item $\ST_{x}(p):=P(x)$
\item $\ST_{x}(\ba):=\ba$
\item $\ST_{x}(\phi \star \psi):=\ST_{x}(\phi) \star \ST_{x}(\psi)$ for $\star\in\{\vee,\wedge,\rightarrow,-\}$
%\item $\ST_{x}(\phi\leq\psi):= \ST_{x}(\phi)\preceq \ST_{x}(\psi)$
%\item $\ST_{x}\Big((\phi_{1}\leq\psi_{1}\& \hdots \&\phi_{n}\leq\psi_{n})\Rightarrow (\phi\leq\psi)\Big):= \Big(\ST_{x}(\phi_{1}\leq\psi_{1})\wedge \hdots \wedge \ST_{x}(\phi_{n}\leq\psi_{n}))\Big) \Rightarrow \\ \ST_{x}(\phi\leq\psi)$
\item $\ST_{x}(\Diamond \psi):=\exists y\left(Rxy \wedge \ST_{y}(\psi)\right)$
\item $\ST_{x}(\Box \psi):=\forall y\left(Rxy \rightarrow \ST_{y}(\psi)\right)$
\item $\ST_{x}(\Db \psi):=\exists y\left(Ryx \wedge \ST_{y}(\psi)\right)$
\item $\ST_{x}(\Bb \psi):=\forall y\left(Ryx \rightarrow \ST_{y}(\psi)\right)$
\item $\ST_{x}(\i):=(\bc_{\i}=x) \wedge \bC_{\i}$
\item $\ST_{x}(\m):=(\bc_{\m}\neq x) \vee \bC_{\m}$
\end{enumerate}
where $y$ is a fresh variable, $p\in \Prop$, $\i\in \Nom$, $\m\in \Cnom$ and $\ba$ is the truth-value constant that corresponds to $a\in A$. The standard translation is extended to inequalities by specifying that:
\[
\ST_{x}(\phi \leq \psi):= \ST_{x}(\phi) \preceq \ST_{x}(\psi).
\]

The standard translation of $\mvml$ into $\mvfol$ is obtained by simply omitting from the above list the clauses proper to $\emvml$.
\end{definition}

\noindent Note that, in the definition above, the variable $x$ will remain unbound and hence does not occur both free and bound in the translation. Moreover,
every two occurrences of quantifiers bind different variables. Therefore, $\ST_{x}(\phi)$ is a clean first-order formula for any modal formula $\phi \in \emvml$.  

%The standard translation of finite sets of inequalities is
% defined similarly to the $2$-valued case: $$\ST_{x}\Big(\{(\phi_{i}\leq\psi_{i})\mid 1\leq i\leq n\}\Big)
%=\{\ST_{x}(\phi_{i})\preceq \ST_{x}(\psi_{i})\mid 1\leq i\leq n\}.$$

Suppose $\modM=(W,R,V)$ is an $\algA$-model for $\mvml$.
Then the {\em corresponding first-order $\algA$-model $\cfoM=(W,R,\{P\}_{p\in\Prop})$ 
of $\mvfol$} is such that
\begin{enumerate}[label=(\roman*)]
\item $W$ and the $\algA$-fuzzy binary relation $R$ are the same as in $\modM$, and
\item if $P$ corresponds to the proposition variable $p$, then $\P$ is
the $\algA$-fuzzy set such that $\P(w)=V(p,w)$ for all $w \in W$.
\end{enumerate}

\noindent To keep the notation light, we will write $\modM$ instead of  $\cfoM$ when this is unlikely to cause confusion.  The notion of the {\em corresponding first-order $\algA$-model of a frame $\modF$} is analogous, but omits the interpretation of unary predicate symbols corresponding to proposition letters and is therefore identical to $\modF$.

In the hybrid logic setting \cite{Hyb}, a nominal $\i$ serves to name a state in
the universe $W$ since any nominal is true at exactly one state.  
With our use of them in the many-valued setting, nominals are used to 
store two pieces of information, namely, the state named by the nominal as
well as its degree of truth at that state. Similar considerations hold for conominals.
This is reflected in the following definition of a corresponding first-order
model for the extended language.
Suppose $\modM$ is an $\algA$-model for $\emvml$.
Then the {\em corresponding first-order $\algA$-model 
\[\cfoM=(W,R,\{P\}_{p\in\Prop},\{c_\i\}_{\i\in\Nom},\{C_\i\}_{\i\in\Nom},
\{c_\m\}_{\m\in\Cnom},\{C_\m\}_{\m\in\Cnom})\]
of $\emvfol$} additionally satisfies:
\begin{enumerate}[label=(\roman*)]
\item For $\i\in\Nom$,  $\bc_i=w$ and 
$C_\i=V(\i,w)$ where $w$ is the state in $W$ such that $V(\i,w)\in J^{\infty}(\algA)$.
\item For $\m\in\Cnom$, $\bc_\m=w$ and 
$C_\m=V(\m,w)$ where $w$ is the state in $W$ such that $V(\m,w)\in M^{\infty}(\algA)$.
\end{enumerate}

\noindent We now show that the standard translation is faithful to the meaning of modal formulas, i.e., that the truth value of modal formula at a state $w$ in an $\algA$-model is the same as the truth value of its standard translation in the corresponding first-order $\algA$-model when $w$ is assigned to the free variable in the translation.

\begin{prop} \label{Proposition:TruthPreserved}
Let $\phi$ be an $\emvml$-formula, $\modM=(W,R,V)$ an $\algA$-model and $w\in W$.  For any assignment $v$ on $\foM$ with $v(x)=w$
\[
V(\phi,w)=\Vert\ST_{x}(\phi)\Vert_{\cfoM,v}.
\]
\end{prop}

\begin{proof}
The proof is by induction on the complexity of $\phi$. The case when $\phi$ is a proposition letter follows immediately from the definitions of $\cfoM$ and the case for truth constants is also immediate. 

For a nominal $\i$, recall that either $w$ is the unique
element such that $V(\i,w)\in J^{\infty}(\algA)$, in which case $V(\i,w)=C_\i$,
or otherwise $V(\i,w)=0$.  In the first case, $\Vert \bc_\i=x\Vert_{\cfoM,v}=1$ since $v(x)=w$.  Thus, 
$V(\i,w)=C_\i=\Vert\bC_\i\Vert_{\cfoM,v}=1\wedge \Vert\bC_\i\Vert_{\cfoM,}
=\Vert\bc_\i=x\Vert_{\cfoM,v}\wedge \Vert\bC_\i\Vert_{\cfoM,v}=
\Vert\left(\bc_\i=x\right)\wedge\bC_\i\Vert_{\cfoM,v}=\Vert \ST_{x}(\i)\Vert_{\cfoM,v}$. In the second case, if $V(\i,w)=0$, then $\Vert\bc_\i=x\Vert_{\cfoM,v}=0$
and we have that $V(\i,w)=\Vert\bc_\i=x\Vert_{\cfoM,v}=\Vert\bc_\i=x\Vert_{\cfoM,v}\wedge \Vert\bC_\i\Vert_{\cfoM,v}
=\Vert\left(\bc_\i=x\right)\wedge\bC_\i\Vert_{\cfoM,v}=\Vert \ST_{x}(\i)\Vert_{\cfoM,v}$.

The case for conominals is similar.

The inductive cases involving $\wedge,\vee,\rightarrow,-$ are straightforward 
and are left to the reader.  We prove only inductive cases involving $\Diamond$
and $\Bb$.  The cases involving $\Box$ and $\Db$ follow similar arguments.

Suppose $\phi=\Diamond \psi$.  Then,
\begin{align*}
\Vert \ST_{x}(\Diamond \psi)\Vert_{\cfoM,v} 
&= \Vert\exists y\left(R xy \wedge \ST_{y}(\psi)\right)\Vert_{\cfoM,v}\\
&=\bigvee \left\{\Vert R xy \wedge \ST_{y}(\psi)\Vert_{\cfoM,v'}\mid  v'\sim_y v\right\}\\
&=\bigvee \left\{R(w,u) \wedge V(\psi,u)\mid u\in W\right\}\\
&=V(\Diamond \psi,w),
\end{align*}
where, to obtain the third equality, we apply the inductive hypothesis and the fact that $x$ is different from $y$ and hence $v'(x) = w$ for all $v'\sim_y v$.

Suppose $\phi=\Bb \psi$.  Then, again applying the inductive hypothesis and fact that $x$ is different from $y$ in the third equality, we have 
\begin{align*}
\Vert\ST_{x}(\Bb \psi)\Vert_{\cfoM,v} 
&= \Vert\forall y\left(R(y,x) \rightarrow \ST_{y}(\psi)\right)\Vert_{\cfoM,v}\\
&=\bigwedge\left\{\Vert R(y,x) \rightarrow \ST_{y}(\psi)\Vert_{\cfoM,v'}\mid v'\sim_y v\right\}\\
&=\bigwedge\left\{R(u,w) \rightarrow V(\psi,u)\mid  u\in W\right\}\\
&=V(\Bb \psi,w).
\end{align*}

\end{proof}

\begin{cor}\label{cor:ineq:faithful:under:ST}
Let $\modM$ be an $\algA$-model, $a\in A$, and $\phi$, $\psi$, $\phi_1, \ldots, \phi_n$ and  $\psi_1, \ldots, \psi_n$ formulas of $\emvml$. Then, for $w\in W$,
\begin{enumerate}
    \item $\modM,w \Vdash_{a} \phi \quad \mbox{ iff } \quad \cfoM \models_{a} \ST_{x}(\phi) [x := w]$,
    \item $\modM \Vdash \phi \leq \psi \quad \mbox{ iff } \quad \cfoM \models \forall x \ST_{x}(\phi \leq \psi)$ and 
    \item $\modM \Vdash \phi_1 \leq \psi_1 \& \cdots \& \phi_n \leq \psi_n  \Rightarrow \phi \leq \psi \quad \mbox{ iff } \quad \cfoM \models \forall x \ST_{x}(\phi_1 \leq \psi_1) \wedge \cdots \wedge \forall x \ST_{x}(\phi_n \leq \psi_n) \to \forall x \ST_{x}(\phi \leq \psi)$. 
\end{enumerate}
\end{cor}

\noindent This has the following consequences for validity on frames: 

\begin{cor}\label{Corollary:StandardTranslation}
Let $\modF$ be an $\algA$-frame, $a\in A$ and $\phi$ a formula of $\emvml$. Then, for $w\in W$,
\begin{enumerate}
    \item $\modF,w \Vdash_{a} \phi \quad \mbox{ iff } \quad \modF \models_{a} \forall \oP\forall \overline{\bc_\i}
	\forall \overline{\bc_\m}\forall \overline{\bC_\i}\forall \overline{\bC_\m}
	\left(\ST_{x}(\phi)\right) [x \coloneqq w]$, 
    \item $\modF \Vdash_{a} \phi \quad \mbox{ iff } \quad \modF \models_{a} \forall \oP\forall \overline{\bc_\i}
	\forall \overline{\bc_\m}\forall \overline{\bC_\i}\forall \overline{\bC_\m}
	\left(\ST_{x}(\phi)\right)$,
    \item $\modF \Vdash \phi \leq \psi \quad \mbox{ iff } \quad \modF \models \forall \oP\forall \overline{\bc_\i}
	\forall \overline{\bc_\m}\forall \overline{\bC_\i}\forall \overline{\bC_\m}
	\left(\forall x \ST_{x}(\phi \leq \psi)\right)$, and
    \item $\modF \Vdash \phi_1 \leq \psi_1 \& \cdots \& \phi_n \leq \psi_n  \Rightarrow \phi \leq \psi \quad \mbox{ iff } \quad\\ \modF \models \forall \oP\forall \overline{\bc_\i}
	\forall \overline{\bc_\m}\forall \overline{\bC_\i}\forall \overline{\bC_\m}
	\left(\forall x \ST_{x}(\phi_1 \leq \psi_1) \wedge \cdots \wedge \forall x \ST_{x}(\phi_n \leq \psi_n) \to \forall x \ST_{x}(\phi \leq \psi)\right)$.
\end{enumerate}
%where $\oP$ is the sequence of predicate variable corresponding to propositional variables occurring in $\phi$,  $\overline{\bc_\i}$ and  $\overline{\bc_\m}$ are, respectively, the sequences of individual variables corresponding to the the nominals and co-nominals occurring in $\phi$, and $\overline{\bC_\i}$ and $\overline{\bC_\m}$ are, respectively, the sequences of predicate variables corresponding to the the nominals and co-nominals occurring in $\phi$.
%
where $\oP$, $\overline{\bc_\i}$ and  $\overline{\bc_\m}$, $\overline{\bC_\i}$ and $\overline{\bC_\m}$ are sequence of, respectively, the  predicate variable corresponding to propositional variables, individual variables corresponding to the nominals, individual variables corresponding to the co-nominals, predicate variables corresponding to the the nominals and predicate variables corresponding to the the co-nominals occurring in the arguments of $\ST_x$.
\end{cor}

\section{Many-valued correspondence theory}\label{sec:mv:corresp:theory}

The classical notion of frame correspondence generalizes to the many-valued setting in the following natural way. Given $a \in A$, an $\mvml$-formula $\phi$ and a $\mvfol$-formula $\alpha(x)$ in one free variable $x$ are {\em local frame $a$-correspondents} if, for all $\algA$-frames $\modF = (W,R)$ and $w \in W$, it holds that
\[
\modF,w \Vdash_{a} \phi \quad \iff \quad \modF \models_{a} \alpha[x:=w].
\label{acorrespondents}
\]
That is, $\phi$ is $a$-valid at $w$ in $\modF$ iff $\alpha$ is $a$-true in $\modF$ under all assignments sending $x$ to $w$. Unpacking the definitions of $a$-validity and $a$-truth, the above is equivalent to requiring that, for all frames $\modF$, it holds that $V(\phi, w)\geq a$ for all valuations $V$ on $\modF$ iff $\Vert \alpha \Vert_{\modF,v} \geq a$ for all assignments $v$ on $\modF$ for which $v(x)=w$.  An $\mvml$-formula $\phi$ and a $\mvfol$-sentence $\alpha$ are {\em global frame $a$-correspondents} if, for all $\algA$-frames $\modF$, it is the case that
\[
\modF\Vdash_{a} \phi \quad \iff \quad \modF \models_{a} \alpha.
\label{Glob:acorrespondents}
\]
Note that it is always the case that $V(\phi,w)\geq 0$ and $v(\alpha)\geq 0$. Hence all many-valued modal formulas are local frame $0$-correspondents to all many-valued first-order formulas, rendering $0$-correspondence trivial. We further note that, via the equivalences in \eqref{eq:satisfac:eqivs}, the notions of local and global $a$-correspondence apply also to inequalities. 

Although there are no general correspondence results in the literature on many-valued modal logic, correspondences have been obtained for a number of particular formulas in \cite{Fitting,Frank}. Most notably, Frankowski \cite{Frank} considers a generalized notion of correspondence which is similar to ours but, instead of demanding that both truth values be at least $\alpha$, he requires that both belong to a given filter $F$ of $\mathbb{A}$.  Restricted to the use of only principle filters --- and therefore to our notion of global correspondence --- Frankowski's generalized conditions are as follows: Let $\modF=(W,R)$ be an $\algA$-frame,  $w_{1},w_{2}\in W$, $a\in A$.  Then
$\mathfrak{F}$ is said to be
\begin{enumerate}[label=(\roman*)]
\item {\em $a$-reflexive} if $\Vert \forall x R(x,x)\Vert_{\modF} \geq a$;
\item {\em $a$-symmetric} if $\Vert\forall x\forall y\left(R(x,y)\rightarrow R(y,x)\right)\Vert_{\modF} \geq a$;
\item {\em $a$-transitive} if 
$\Vert \forall x\forall y\forall z\left(\left(R(x,y)\wedge R(y,z)\right)\rightarrow R(x,z)\right)\Vert_{\modF}\geq a$;
\item {\em $a$-dense} if 
$\Vert\forall x\forall y\left(R(x,y)\rightarrow \exists z\left(R(x,z) \wedge R(z,y)\right)\right)\Vert_{\modF} \geq a$;
\item {\em $a$-right unbounded} (or {\em $a$-serial}) if $\Vert \forall x\exists y R(x,y)\Vert_{\modF} \geq a$.
\end{enumerate}

It is shown \cite{Frank} that the well-known modal formulas which correspond, respectively, to reflexivity, symmetry, transitivity, density and right unboundedness in the 2-valued case also correspond, respectively, to the the $a$-versions of these properties in the many-valued case.  Is this part of a more general phenomenon? Specifically, in our study of many-valued correspondence theory, we aim to answer the following questions:
\begin{enumerate}[label=(\roman*)]
\item Can a modal formula $\phi$ in $\mvml$ and a first-order formula $\alpha$ in $\mvfol$ be local 
frame $a_1$-correspondent, but not local $a_2$-correspondents for $a_1,a_2\in A$ with $a_1\neq a_2$? A related question is:  can a modal formula have distinct correspondents depending on the truth-value?

\item Can we identify a syntactically described class of modal formulas with (effectively computable) many-valued local frame correspondents in $\mvfol$?

\item If a modal formula has a many-valued frame correspondent, how does this correspondent compare to
it correspondent in the two-valued setting?
\end{enumerate}

\noindent The first of these questions is answered in the affirmative in the example below.  That is, we illustrate
how one modal formula can have different many-valued correspondents, depending on the truth-values.
The answer to the second question is the focus of Sections~\ref{sec:Sahl:and:Ind:fmls} and \ref{Section:ALBA}, where, in Section~\ref{sec:Sahl:and:Ind:fmls}, we define the classes of inductive \cite{Goranko2006-GORECF} and Sahlqvist formulas for $\mvml$ via the methodology developed in \cite{ALBA, nonDistALBA} and, in Section \ref{Section:ALBA}, introduce a many-valued version the ALBA-algorithm \cite{ALBA, nonDistALBA} which we use to prove that all inductive (and hence all Sahlqvist) $\mvml$-formulas  have many-valued first-order correspondents. Finally, in Section~\ref{Section:SvBResult}, we show that many-valued correspondent of any classical Salhqvist formula is syntactically identical to its classical correspondent.

\begin{figure}[h]%
	\begin{center}
	\begin{minipage}[c]{0.3\textwidth}
	\begin{tikzpicture}[node/.style={circle, draw, fill=black}, scale=0.8]
		\node[circle,fill=black!100,minimum size=5pt,inner sep=0pt, label=below:$0$] (0)  at (0, -2) {};
		\node[circle,fill=black!100,minimum size=5pt,inner sep=0pt, label=left:$\alpha$] (alpha)  at (-1, -1) {};
		\node[circle,fill=black!100,minimum size=5pt,inner sep=0pt, label=right:$\beta$] (beta)  at (1, -1) {};	
		\node[circle,fill=black!100,minimum size=5pt,inner sep=0pt, label=right:$\gamma$] (gamma)  at (0, 0) {};	
		\node[circle,fill=black!100,minimum size=5pt,inner sep=0pt, label=above:$1$] (1)  at (0, 1) {};	
		%Edges
		\draw[-] (0) to node {} (alpha);
		\draw[-] (0) to node {} (beta);
		\draw[-] (alpha) to node {} (gamma);
		\draw[-] (beta) to node {} (gamma);
		\draw[-] (gamma) to node {} (1);
	\end{tikzpicture}
	\end{minipage}
	%\hspace{2cm}
	\begin{minipage}[c]{0.3\textwidth}
	\begin{tabular}{r|ccccc}
		$\to$ &$0$ &$\alpha$ &$\beta$ &$\gamma$ &$1$\\
		\hline
		$0$ &$1$ &$1$ &$1$ &$1$ &$1$\\
		$\alpha$ &$\beta$ &$1$ &$\beta$ &$1$ &$1$\\
		$\beta$ &$\alpha$ &$\alpha$ &$1$ &$1$ &$1$\\
		$\gamma$ &$0$ &$\alpha$ &$\beta$ &$1$ &$1$\\
		$1$ &$0$ &$\alpha$ &$\beta$ &$\gamma$ &$1$
	\end{tabular}
	\end{minipage}
	%\hspace{2cm}
	\begin{minipage}[c]{0.3\textwidth}
	\begin{tabular}{r|l}
		$a$ &$\neg a$\\
		\hline
		$0$ &$1$\\
		$\alpha$ &$\beta$\\
		$\beta$ &$\alpha$\\ 
		$\gamma$ &$0$\\
		$1$ &$0$\
	\end{tabular}
	\end{minipage}
	\end{center}
	\caption{The Heyting algebra $\mathbb{P}$ with tables for its implication and negation.}\label{HAP}
\end{figure}
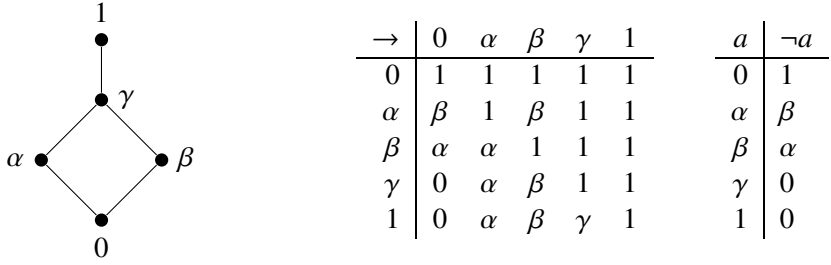

\begin{example}\label{example:change:from:classical:corresp}
Let $\mathbb{P}=(\{0,\alpha,\beta,\gamma,1\},\vee,\wedge,\rightarrow,0,1)$ be
the Heyting algebra depicted in Figure~\ref{HAP}. Let $\modF=(W,R)$ be a $\mathbb{P}$-frame.  We show first that
$\neg p \vee \Diamond p$ and $\mathbf{0}$ are local frame $1$-correspondents. Then, in contrast to this, we show that
$\neg p \vee \Diamond p$ and $a$-reflexivity are local frame $a$-correspondents, for $a\neq 1$, i.e., that $\modF,w \Vdash_{a} \neg p \vee \Diamond p$ iff $\modF \models_{a} Rxx [x:=w]$.

\begin{enumerate}[label=(\roman*)]
\item Let $w\in W$ and suppose $V$ is a valuation on $\modF$ such that
$V(p,v)=\alpha$ for all $v\in W$.  Then,
\[
V(\neg p \vee \Diamond p,w)= V(\neg p,w) \vee \bigvee\left\{R(w,v) \wedge V(p,v)\mid v\in W\right\}
\leq \beta \vee \alpha= \gamma.
\]
So $\mathfrak{F},w \not\Vdash_{1} \neg p\vee \Diamond p$. Also, $\Vert \mathbf{0} \Vert_{\modF} = 0 \not\geq 1$, so $\mathfrak{F} \not\models_{1} \mathbf{0} [x:= w]$. Hence, $\neg p \vee \Diamond p$ and $\mathbf{0}$ are local frame $1$-correspondents.

\item We prove the second claim for $a=\gamma$.
The proof is similar if $a=\alpha$ or $a=\beta$.
\begin{enumerate}[label=(\alph*)]
\item Suppose first that $R(w,w)\ngeq \gamma$, i.e., $w$ is not $\gamma$-reflexive.
Now let $V$ be a valuation on $\modF$ such that $V(p,w)=1$ and $V(p,u)=0$ for all $u\in W$ such that $u\neq w$.
Then,
\[
V(\neg p \vee \Diamond p,w)= V(\neg p,w) \vee \bigvee\left\{R(w,u) \wedge V(p,u)\mid u\in W\right\}
< 0 \vee \gamma.
\]
The inequality follows since the maximum value of $\bigvee\left\{R(w,u) \wedge V(p,u)\mid u\in W\right\}$
is $\alpha \wedge 1$ when $R(w,w)=\alpha$ or $\beta \wedge 1$ when $R(w,w)=\beta$.
Hence, $\modF,w \not\Vdash_{\gamma} \neg p\vee \Diamond p$.

\item Next suppose that $R(w,w)\geq \gamma$, i.e. $w$ is $\gamma$-reflexive.
Let $V$ be any valuation on $\modF$. Then,
\begin{align*}
V(\neg p \vee \Diamond p,w)&= V(\neg p,w) \vee \bigvee\left\{R(w,u) \wedge V(p,u)\mid u\in W\right\}\\
&\geq V(\neg p,w) \vee \left(\gamma \wedge V(p,w)\right)\\
&= \left(V(\neg p,w) \vee \gamma\right) \wedge \left(V(\neg p,w) \vee V(p,w)\right)\\
&\geq \gamma.
\end{align*}
The final inequality follows from the fact that $\neg a \vee a \geq \gamma$ for all $a \in \mathbb{P}$. 
Hence, $\modF,w \Vdash_{\gamma} \neg p\vee \Diamond p$.
\end{enumerate}  
\end{enumerate}
\end{example}

\section{Sahlqvist and inductive formulas}\label{sec:Sahl:and:Ind:fmls}

A line of research known as \emph{unified correspondence} (see e.g.\ \cite{ALBA,CoGhPa14,nonDistALBA}), proposes a general definition for classes of inductive and Sahlqvist formulas that uniformly applies to all logics, called \emph{LE-logics}, algebraically captured by classes of normal lattice expansions. This definition is based purely on the order-theoretic properties of the algebraic operations interpreting the logical connectives. All such formulas are canonical and, moreover, define first-order properties on any relational semantics connected to the algebraic semantics via a suitable discrete duality. To be able to apply these definitions in the present context, it is enough to note that the complex algebras of $\algA$-frames, being Heyting algebras with normal operators (cf.\ Lemma \ref{lem:complex:is:perf:Heyting}) are, in particular, normal lattice expansions.

%\subsection{The Sahlqvist and inductive formulas: a general purpose definition}\label{subsec:SahlqvistGeneralPurpose}

The definition of Sahlqvist and inductive formulas for logics lacking classical negation in necessarily more involved than for logics with such a negations. In particular, where definitions in the latter setting may, without loss of generality, privilege e.g.\ positive occurrences of variables in `constructive' definitions, in the former setting careful bookkeeping of polarities and `forbidden path' definitions are required.  This, and tracking the nesting of connectives, is  accomplished using the notions signed generation trees and polarity types, first used for similar purposes in \cite{GEHRKE200565}.     

%\texttt{Write an introductory paragraph. Lack of classical negation necessitates bookkeeping about polarities, hence generation trees (cite GNV)....}\wnote{Willem to complete this section. We are still lacking the definition of the inductive formulas (good branches and all we need for that) and examples.}

Based on the usual notion of a \emph{generation tree} of a formula, a {\em signed generation tree} (see, e.g., \cite{GEHRKE200565}) further associates a sign ($+$ or $-$) with each node. We assign a sign to the root of the generation tree and this is then propagated to the other nodes in the tree as follows: children of nodes labelled with all connectives \emph{except} $\rightarrow$ receive the same sign as their parent, while the left child of a node labelled $\rightarrow$ receives the opposite to its parent's while the right child inherits the same sign as its parent.   The \emph{positive (negative) generation tree of} a formulas $\phi$, denoted $+\phi$ ($- \phi$), is obtained by signing the of $\phi$'s generation tree with $+$ ($-$) and propagating the signs in the manner described.

An occurrence of a propositional variable $p$ is \emph{positive} (\emph{negative}) \emph{in a formula $\phi$} if the corresponding node is signed $+$ ($-$) in the positive generation tree of $\phi$. The formula $\phi$ \emph{is positive (negative) in $p$} if all occurrences of $p$ in $\phi$ are positive (negative).

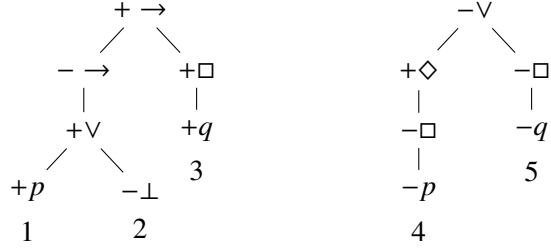
\begin{figure}
\begin{center}
\begin{tikzpicture}[level distance=8mm]
%\tikzset{style={align=center,anchor=north}}
\tikzstyle{level 6}=[sibling distance=30mm]
\node{$+ \to$}
        child{node {$-\to$}
            child{node[label=below:1] {$+ p$}}
            child{node[label=below:2] {$- \bot$}}
            }
        child{node {$+ \Box$}
            child{node[label=below:3]  {$+q$}}
        };
\end{tikzpicture}
\hspace{2cm}
\begin{tikzpicture}[level distance=8mm]
%\tikzset{style={align=center,anchor=north}}
\tikzstyle{level 6}=[sibling distance=30mm]
\node {$-\vee$}
            child{node {$-\Diamond$}
                child{node {$- \Box$}
                    child{node[label=below:4] {$-q$} }
                }
            }
        child{node {$-\Box$} 
            child{node[label=below:5] {$-p$} }
            };
\end{tikzpicture}
\end{center}
\caption{The generation trees of $+ ((p \to \bot) \to \Box q)$ and $- (\Diamond \Box q \vee \Box p)$ with branches numbered for ease of reference.}\label{Fig:Sahl:fml}

\end{figure}

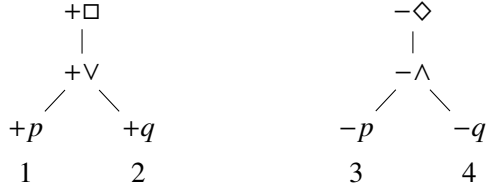
\begin{figure}
\begin{center}
\begin{tikzpicture}[level distance=8mm]
%\tikzset{style={align=center,anchor=north}}
\tikzstyle{level 6}=[sibling distance=30mm]
\node{$+ \Box$}
      child{node {$+\vee$}
            child{node[label=below:1] {$+p$}}
            child{node[label=below:2] {$+q$}}
            };
\end{tikzpicture}
\hspace{2cm}
\begin{tikzpicture}[level distance=8mm]
%\tikzset{style={align=center,anchor=north}}
\tikzstyle{level 6}=[sibling distance=30mm]
\node {$-\Diamond$}
            child{node {$-\wedge$}
                child{node[label=below:3] {$-p$}}
                child{node[label=below:4] {$-q$}}
                };
\end{tikzpicture}
\end{center}
\caption{The generation trees of $+ \Box(p \vee q)$ and $- \Diamond(p \wedge q)$ with branches numbered for ease of reference.}\label{Fig:ind:non:Sahl}
\end{figure}

\begin{figure}
\begin{center}
\begin{tikzpicture}[level distance=8mm]
%\tikzset{style={align=center,anchor=north}}
\tikzstyle{level 1}=[sibling distance=25mm]
\tikzstyle{level 2}=[sibling distance=15mm]
\node{$- \to$}
      child{node {$+\wedge$}
            child{node {$+\Box$}
                child{node {$+ \to$}
                    child{node {$- \wedge $}
                        child{node[label=below:1] {$-{\ba}$}}
                        child{node[label=below:2] {$-p$}}
                        }
                    child{node[label=below:3] {$+q$}}
                    }
                }
            child{node {$+\Box$}
                child{node[label=below:4] {$+p$}}
                }
            }
        child{node {$-\Diamond$}
            child{node {$-\Box$}
                child{node[label=below:5] {$-q$}}
                }
            }    
            ;
\end{tikzpicture}
\end{center}
\caption{The negative generation trees of $\Box(\ba \wedge p \to q) \wedge \Box p \to \Diamond \Box q$ with branches numbered for ease of reference.}\label{Fig:induct:with:Truth:Constant}
\end{figure}
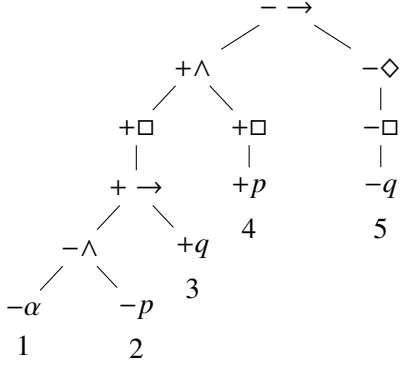

\begin{example}\label{ex:signed:gen:trees}
The positive generation tree of $(p \to \bot) \to \Box q$ and the negative generation tree of $\Diamond \Box q \vee \Box p$  are displayed in Figure \ref{Fig:Sahl:fml}, while Figure \ref{Fig:ind:non:Sahl} displays the positive generation tree of $\Box(p \vee q)$ and negative generation tree of $ \Diamond(p \wedge q)$. Lastly, the negative generation tree of the formula $\Box(\bm{\alpha} \wedge p \to q) \wedge \Box p \to \Diamond \Box q$ for $\alpha \in \algA$ is given in Figure \ref{Fig:induct:with:Truth:Constant}.
\end{example}

\begin{definition}[Order types and critical branches]
    An {\em order type} over $n\in \mathbb{N}$ is an $n$-tuple $\epsilon\in \{1, \partial\}^n$. We will often associate  order types with tuples (or sets) of variables $(p_1,\ldots p_n)$. In this case we may think of $\epsilon$ as a map $\epsilon : \{ p_1,\ldots p_n\} \to \{1, \partial\}^n$ and accordingly write $\epsilon(p_i)$ instead of $\epsilon_i$.
    
    Given an $\mvml$-formula $\phi(p_1,\ldots p_n)$, an order type $\epsilon$ over $n$, and $1 \leq i \leq n$, a leaf in a signed generation tree of $\phi$ labelled $+p_i$ with $\epsilon_i = 1$ or $-p_i$ with $\epsilon_i = \partial$ is an \emph{$\epsilon$-critical node}. A branch ending in an $\epsilon$-critical node is an $\epsilon$-{\em critical branch}.
\end{definition}
\noindent We are now ready to start defining the classes of Sahlqvist and inductive formulas in the setting of many-valued modal logic. These definitions are instantiations for many-valued modal logic of those for generic LE-logics found in \cite{nonDistALBA}.  

\begin{table}[h]
	\begin{center}
		\bgroup
		\def\arraystretch{1.5}
		\begin{tabular}{| c | c |}
			\hline
			Skeleton-nodes &PIA-nodes\\
			\hline\hline
			$\Delta$-adjoints & Syntactically Right Adjoints (SRA)\\
			\hline
			\bgroup
			\def\arraystretch{1.5}
			\begin{tabular}{ c c c c c c c c  }
				$+$ & $\quad$ & $\vee$ &$\wedge$ &${}$ &${}$  &${}$ &${}$\\
				$-$ & $\quad$ & $\wedge$ &$\vee$\\
			\end{tabular}
			\egroup
			&\bgroup
			\def\arraystretch{1.5}
			\begin{tabular}{ccccc}
				$+$ & $\quad$&$\wedge$ &$\Box$ &${}$\\
				$-$ & $\quad$&$\vee$ &$\Diamond$ &${}$\\
			\end{tabular}
			\egroup
			\\
			\hline
			Syntactically Left Residuals (SLR) & Syntactically Right Residuals (SRR)\\
			\hline
			\bgroup
			\def\arraystretch{1.5}
			\begin{tabular}{ c c c c c c c c  }
				$+$ &$\quad$ &$\Diamond$  &${}$ &${}$  &${}$ &${}$  &${}$\\
				$-$ & $\quad$&$\Box$ &$\rightarrow$\\
			\end{tabular}
			\egroup
			&\bgroup
			\def\arraystretch{1.5}
			\begin{tabular}{cc c c c}
				$+$ & $\quad$&$\vee$  &$\rightarrow$\\
				$-$ & $\quad$&$\wedge$ \\
			\end{tabular}
			\egroup\\
			\hline
		\end{tabular}
		\egroup
	\end{center}
	\caption{Classification of nodes}\label{Join:and:Meet:Friendly:Table}
\end{table}

\begin{definition}[Node classification, good and excellent branches, cf.\ {\cite[Definition 3.2]{nonDistALBA}}]\label{Excellent:Branch:Def}
    Nodes in signed generation trees of $\mvml$-formulas are classified according to Table \ref{Join:and:Meet:Friendly:Table} into \emph{skeleton} and \emph{PIA}-nodes\footnote{The acronym PIA stands for `positive implies atomic' and was introduced by van Benthem in in \cite{van2005minimal}. Intuitively, skeleton nodes form the `framework' or outer part of formulas which hold the inner, PIA parts.} which are, respectively, further subdivided into \emph{$\Delta$-adjoints} and \emph{syntactically left residuals (SLR)}, and \emph{syntactically right adjoints (SRA)} and \emph{syntactically right residuals (SRR)}.  
    
    A branch (i.e., a path from the root to a leaf) in a signed generation tree is \emph{good} if it is the concatenation of two paths, $P_1$ and $P_2$, one of which may possibly be of length $0$, such that $P_1$ is a path from the root\footnote{In some sources, e.g. \cite{nonDistALBA}, the direction of the paths is the other way around, but the present definition is equivalent to that of \cite{nonDistALBA} when instantiated to $\mvml$.} consisting (apart from a possible variable or constant node) only of skeleton-nodes, and $P_2$ consists (apart from a possible variable or constant node) only of PIA-nodes. A good branch is \emph{excellent} if $P_2$ consists only of SRA-nodes.
\end{definition}
\begin{remark}
    Note that some nodes are both Skeleton and PIA. The transition, on a good branch, from its $P_1$ to its $P_2$ part may therefore not be uniquely determined. However, on such a branch, the first occurrence of a PIA-only node (i.e., of $+\Box$, $-\Diamond$ or $+ \to$) indicates the latest possible  transition from $P_1$ to $P_2$.
\end{remark}

\begin{example}\label{ex:classfctn:branches}
    As the reader can easily verify, branches 1, 2 and 4 in Figure \ref{Fig:Sahl:fml} are not good, while branches 3 and 5 are excellent.  In Figure \ref{Fig:ind:non:Sahl} all four are  good but not excellent. Lastly, if Figure \ref{Fig:induct:with:Truth:Constant} branches 1, 2 and 3 are good (but not excellent), branch 4 is excellent and branch 5 is not good. 
\end{example}

\noindent We can now define the Sahlqvist inequalities and formulas. The definition of the more general class of inductive inequalities (and formulas) is more involved and will be given after. 
\begin{definition}[Shalqvist inequalities, cf.\ {\cite[Definition 3.5]{nonDistALBA}}]\label{Sahlqvist:Ineq:Def}
	For any order type $\epsilon$, the signed generation tree of a formula $\phi$ is \emph{$\epsilon$-Sahlqvist} if every $\epsilon$-critical branch is excellent. An inequality $\phi \leq \psi$ is \emph{$\epsilon$-Sahlqvist} if the trees $+ \phi$ and $- \psi$ are both $\epsilon$-Sahlqvist.  An inequality is \emph{Sahlqvist} if it is $\epsilon$-Sahlqvist for some $\epsilon$. %A formula $\psi$ is \emph{Sahlqvist} if if the inequality $\top \leq \psi$ is Sahlqvist.
\end{definition}

\begin{definition}[Sahlqvist formulas]\label{Sahlqvist:Form:Def}
	A formula $\psi$ is a \emph{Sahlqvist formula} if $\top \leq \psi$ is a Sahlqvist inequality.
\end{definition}

\begin{example}\label{ex:Sahl:and:non:Sahl}
    Consider the inequality $((p \to \bot) \to \Box q)  \leq \Diamond\Box q \vee \Box p$. The positive generation tree of the left hand side and negative generation tree of the right hand side were given in Figure \ref{Fig:Sahl:fml} and we noted in Example \ref{ex:classfctn:branches} that the branches ending in leaves $+p$, $-\bot$ and $-q$ we not good while those ending in $+q$ and $-p$ we excellent. The inequality is therefore $\epsilon$-Sahlqvist for $\epsilon(p) = \partial$ and $\epsilon(q) = 1$. 

    Next consider the inequality $\Box(p \vee q) \leq \Diamond (p \wedge q)$. The positive and negative generation trees of its left hand right hand sides, respectively, were displayed in Figure  \ref{Fig:ind:non:Sahl}.  In Example \ref{ex:classfctn:branches} we noted that all the branches in these trees were good but not excellent. Therefore this inequality is not Sahlqvist as it fails to be $\epsilon$-Sahlqvist for any order type $\epsilon$.  

    The negative generation tree of the formula $\Box(\ba \wedge p \to q) \wedge \Box p \to \Diamond \Box q$ was given in Figure \ref{Fig:induct:with:Truth:Constant}. As indicated in Example \ref{ex:classfctn:branches}, the branch ending in $+q$ is good but not excellent and the branch ending in $-q$ is not good. There is therefore now way of choosing an order type $\epsilon$ such that all $\epsilon$-critical branches are excellent. Therefore the inequality $\top \leq \Box(\ba \wedge p \to q) \wedge \Box p \to \Diamond \Box q$ and hence the formula $\Box(\ba \wedge p \to q) \wedge \Box p \to \Diamond \Box q$ are not Sahlqvist.
\end{example}

\begin{definition}[Inductive inequalities, cf.\ {\cite[Definition 3.4]{nonDistALBA}}]\label{Inducive:Ineq:Def}
Given an order type $\epsilon$ and a strict partial order $<_{\Omega}$ on the variables $p_1,\ldots p_n$, the signed generation tree $\ast \phi$, $\ast \in \{-, + \}$, of a formula $\phi(p_1,\ldots p_n)$ is \emph{$(\Omega, \epsilon)$-inductive}
if, for all $1 \leq i \leq n$, every $\epsilon$-critical branch with leaf labelled $p_i$ is good and, moreover,  for every subtree $\star (\chi \circledast \zeta)$ rooted at an SRR node $\star \circledast  \in \{ +\vee, +\to,  -\wedge \}$ on such a branch, either
\begin{enumerate}
\item there are no critical branches through the signed subtree for $\chi$ and $p_j <_{\Omega} p_i$ for every $p_j$ occurring in $\chi$, or 
\item there are no critical branches through the signed subtree for $\zeta$ and $p_j <_{\Omega} p_i$ for every $p_j$ occurring in $\zeta$.
\end{enumerate}
We will refer to $<_{\Omega}$ as  the \emph{dependency order} on the variables. An inequality $\phi \leq \psi$ is \emph{$(\Omega, \epsilon)$-inductive} if the trees $+\phi$ and $-\psi$ are both $(\Omega, \epsilon)$-inductive.  An inequality $\phi \leq \psi$ is \emph{inductive} if it is $(\Omega, \epsilon)$-inductive for some $\Omega$ and $\epsilon$. A formula $\phi$ is inductive when the inequality $\top \leq \phi$ is or, equivalently, when the tree $-\phi$ is $(\Omega, \epsilon)$-inductive for some $\Omega$ and $\epsilon$.
\end{definition}

Note that since excellent branches are good and do not contain any SRR-nodes, it follows from the definitions that all Sahlqvist inequalities/formulas are inductive.  

\begin{example}
    In Example \ref{ex:Sahl:and:non:Sahl} we noted that the inequality  $\Box(p \vee q) \leq \Diamond (p \wedge q)$ was not Sahlqvist. We argue that it is also not inductive. Indeed, suppose that it were $(\Omega, \epsilon)$-inductive for some $\Omega$ and $\epsilon$. Since $+p$ and $+q$ are leaves on branches though the same SRR-node $+\vee$
    \footnote{Note that we have to classify $+\vee$ as an SRR-node for the two branches passing through it to be good.} we will have to have either (a) $\epsilon(p) = 1$ and $\epsilon(q) = \partial$ with $q <_{\Omega} p$, or (b) $\epsilon(p) = \partial$ and $\epsilon(q) = 1$ with $p <_{\Omega} q$. But, given that there are also leaves $-p$ and $-q$ on branches through the same SRR-node $-\wedge$, having $\epsilon(q) = \partial$ would require $p <_{\Omega} q$ making (a) nonviable or, having  $\epsilon(p) = \partial$ would require $q <_{\Omega} p$ making (b) nonviable. Thus no such pair $(\Omega, \epsilon)$ can exist.

    Now consider the formula $\Box(\ba \wedge p \to q) \wedge \Box p \to \Diamond \Box q$, shown in Example \ref{ex:Sahl:and:non:Sahl} to be non-Sahlqvist. However, if we take $\epsilon(p) = \epsilon(q) = 1$ and $p <_{\Omega} q$, then the negative generation tree of this formula (Figure \ref{Fig:induct:with:Truth:Constant}) is $(\Omega, \epsilon)$-inductive, as this ensures that the $\epsilon$-critical branches (branches 3 and 4) are good and the the SRR-node $+ \to$ occurring on branch 3 satisfies condition 1 of Definition \ref{Inducive:Ineq:Def}. 
\end{example}

\section{ALBA for many-valued modal logics}\label{Section:ALBA}

In this section we present a version of the ALBA algorithm/calculus \cite{ALBA}, suitably instantiated for the connectives of the language $\mvml$. Apart from small but crucial modifications to the input step and  the translation stage, and the presence of truth value constants in the language, this is, on the syntactic level, the exact same algorithm/calculus as presented in \cite{ALBA} for the language of intuitionistic modal logic. However, the semantic interpretation of the rules is of course different, namely, instead of being interpreted on intuitionistic Kripke frames or Heyting algebras with operators, they are now being interpreted on $\algA$-frames or, dually, on their complex algebras. We show that the rules are sound with respect to this new interpretation and that, therefore, the first-order formulas derived are indeed local frame $a$-correspondents of the input formulas. We can then appeal to the purely syntactic fact that ALBA successfully reduces all inductive formulas (and hence also all Sahlqvist formulas) to conclude that these formulas all have effectively computable first-order $a$-correspondents over $\algA$-frames.

\subsection{The algorithm}

The ALBA algorithm for $\mvml$ tries to derive a frame $a$-correspondent for an $\mvml$-formula or inequality. Given an inequality $\phi \leq \psi$,  we pass to the algorithm the inequality $\phi \wedge \ba \leq \psi$. Or given a formula $\psi$, we pass the inequality $1 \wedge \ba \leq \psi$. On the syntactic level, this addition of the conjunct $\ba$ is really the only difference between the algorithm we present here that given in \cite{ALBA} for the language of intuitionistic modal logic.

\paragraph{Phase 1: Preprocessing and first approximation.} ALBA receives $\phi\wedge\ba \leq\psi$ as input and performs the following preprocessing steps exhaustively on it:
\begin{enumerate}
\item in the left hand side of the inequality, `bubble up' disjunctions as far as possible by applying the equivalences $\Diamond(\alpha \vee \beta) \equiv \Diamond\alpha \vee \Diamond \beta$ and $\gamma \wedge (\alpha \vee \beta) \equiv (\gamma \wedge \alpha) \vee (\gamma \wedge \beta)$;

\item in the right hand side of the inequality, `bubble up' conjunctions as far as possible by applying the equivalences $\Box(\alpha \wedge \beta) \equiv \Box \alpha \vee \Box \beta$,  and $\gamma \vee (\alpha \wedge \beta) \equiv (\gamma \vee \alpha) \wedge (\gamma \vee \beta)$, $(\alpha \vee \beta) \rightarrow \gamma \equiv (\alpha \rightarrow \gamma) \wedge (\beta \rightarrow \gamma)$
\item apply the following splitting rules:
\[\frac{\alpha \leq \beta \wedge \gamma}{\alpha \leq \beta \quad \alpha \leq \gamma} \quad \frac{\alpha \vee \beta \leq \gamma}{\alpha \leq \gamma \quad \beta \leq \gamma};\]
\item apply the elimination of monotone and antitone variable rules, 
\[\frac{\alpha(p) \leq \beta(p)}{\alpha(\mathbf{0})\leq \beta(\mathbf{0})} \quad \frac{\gamma(p)\leq \delta(p)}{\gamma(\mathbf{1})\leq \delta(\mathbf{1})}\]
where $\alpha(p) \leq \beta(p)$ is positive and $\gamma(p) \leq \delta(p)$
is negative in $p$, respectively.
\end{enumerate}
%
%\wmnote{Maybe reformulate rules to work on quasi-inequalities as in RAMICS paper, then we do not need all the side conditions.}
Denote the finite set $\{\phi_{i} \wedge \mathbf{a}\leq \psi_{i} \mid 1 \leq i \leq \ell \}$
of inequalities obtained after the exhaustive application of these rules by 
$\Pre(\phi  \wedge \ba\leq \psi)$.  

Now continue with each inequality $\phi_{i} \wedge \ba\leq \psi_{i}$ in 
$\Pre(\phi  \wedge \ba\leq \psi)$ individually as follows:
\begin{enumerate}
%\wmnote{The examiners asked why only once - perhaps here or later we should explain why only once.}
\item apply the following {\em first-approximation} rule once \[\frac{\phi_{i} \wedge \ba\leq \psi_{i}}
{\i_{0} \leq \phi_{i} \wedge \mathbf{a}\quad \psi_{i} \leq \m_{0}},\]
where $\i_{0}$ is a nominal and $\m_{0}$ is a co-nominal, to obtain a system of inequalities
of the form $\{\i_{0} \leq \phi_{i} \wedge \ba, \psi_{i} \leq \m_{0}\}$
called an {\em initial system};
\item apply the first splitting rule to the inequalities in the initial system obtained in the previous step
to produce a system of the form $\{\i_{0} \leq \phi_{i},\i_{0}\leq \ba, \psi_{i} \leq \m_{0}\}$.
\end{enumerate}
ALBA continues on each system of the form 
$\{\i_{0} \leq \phi_{i},\i_{0}\leq \ba, \psi_{i} \leq \m_{0}\}$ separately.

\paragraph{Phase 2: Reduction and Elimination.} This phase proceeds exactly as in \cite{ALBA}, whereby Residuation
and Approximation rules are applied to the system of inequalities
given by $\{\i_{0} \leq \phi_{i}~,~\i_{0}\leq \mathbf{a}~,~ \psi_{i} \leq \m_{0}\}$.
Such rules will only be applicable to the inequalities $\i_{0} \leq \phi_{i}$ and $\psi_{i} \leq \m_{0}$, leaving $\i_{0}\leq \mathbf{a}$ untouched.  
In this phase, we attempt to remove all propositional variables from a given initial system by applying the Ackermann rules given below.  In order to apply an Ackermann rule, the system must have a specific shape, which can be obtained by applying the following rules, if possible.
	\begin{enumerate}
		\item Residuation rules:  %All operations in the signature of the complex algebra $\mathfrak{F}^{+}$ are either left adjoints or right adjoints of each other.  Hence we have the following \emph{Residuation rules}. 
        \[
        \frac{\phi \wedge \psi\leq \gamma}{\phi \leq \psi\rightarrow \gamma} \quad \frac{\phi \leq \psi \vee \gamma}{\phi - \psi\leq \gamma} \quad \frac{\Diamond \phi\leq \psi}{\phi \leq \Bb \psi} \quad \frac{\phi \leq \Box \psi}{\Db \phi \leq \psi}
        \]
		\item Approximation rules: %The complex algebra $\mathfrak{F}^{+}$ is a perfect Heyting algebra with operators.  The \emph{Approximation rules} follow from the fact that each element is the join of completely join-irreducible elements below it and the meet of completely meet-irreducible elements above it and the operators have infinitary distribution properties. 
		
		\[
        \frac{\i \leq \Diamond \phi}{\j \leq \phi \quad \i \leq \Diamond \j} \quad \frac{\Box \phi \leq \m}{\phi \leq \n \quad \Box \n \leq \m}
        \]
        \[
        \frac{\phi\rightarrow \psi \leq \m}{\j \leq \phi \quad \j \rightarrow \psi \leq \m}\quad \frac{\phi\rightarrow \psi \leq \m}{\psi \leq \n \quad \phi\rightarrow \n \leq \m}
        \]
        where the nominals and co-nominals introduced in the conclusions of the approximation rules must not occur anywhere in the system before the rule is applied.
    \item Ackermann rules: Once application of the Residuation and Approximation rules
     results in a system of the desired shape, the Ackerman rules can be applied to eliminate propositional letters from the system.  The Ackermann rules follow from the respective Ackermann lemmas \cite{ALBA}. In the formulation of the Ackermann rules below, the inequality $\i_0 \leq \mathbf{a}$ need not have been listed separately, since it satisfies the requirements on the $\beta_i \leq \gamma_i$ formulas and could hence been absorbed into that list. However, we keep is separate to emphasize the fact that it is being carried along unchanged.  The soundness of the Ackermann rules follow from the Ackermann lemmas (see \cite{Ackermann:Untersuchung, ALBA}). 
		\begin{enumerate}
		\item The Right-hand Ackermann rule:  The system 
        \[
        \i_{0}\leq \mathbf{a}, \: \alpha_{1}\leq p, \ldots, \alpha_{n}\leq p, \: \beta_{1}\leq \gamma_{1}, \ldots, \beta_{m}\leq \gamma_{m} 
        \]
        is replaced with the system
        \[
        \begin{tabular}{l}
        $\i_{0}\leq \mathbf{a}$,\\ 
        $\beta_{1}((\alpha_1 \vee \cdots \vee \alpha_n)/p)\leq \gamma_{1}((\alpha_1 \vee \cdots \vee \alpha_n)/p)$,\\
        $\vdots$\\
        $\beta_{m}((\alpha_1 \vee \cdots \vee \alpha_n)/p)\leq \gamma_{m}((\alpha_1 \vee \cdots \vee \alpha_n)/p)$ 
        \end{tabular}
        \]
  	where $p$ does not occur in      $\alpha_{1},\hdots,\alpha_{n}$, each $\beta_{i}$ is positive, and each $\gamma_{i}$ is negative in $p$.
    \item The Left-hand Ackermann rule: The system  
		\[
        \i_{0}\leq \mathbf{a}, p\leq \alpha_{1}, \ldots,  p\leq \alpha_{n}, \: \beta_{1}\leq \gamma_{1}, \ldots,\beta_{m}\leq \gamma_{m} 
        \]
         is replaced with the system
        \[
        \begin{tabular}{l}
        $\i_{0}\leq \mathbf{a}$,\\ 
        $\beta_{1}((\alpha_1 \wedge \cdots \wedge \alpha_n)/p)\leq \gamma_{1}((\alpha_1 \wedge \cdots \wedge \alpha_n)/p)$,\\
        $\vdots$\\
        $\beta_{m}((\alpha_1 \wedge \cdots \wedge \alpha_n)/p)\leq \gamma_{m}((\alpha_1 \wedge \cdots \wedge \alpha_n)/p)$ 
        \end{tabular}
        \]    
        where 
        $p$ does not occur in $\alpha_{1},\hdots,\alpha_{n}$, each $\beta_{i}$ is negative in $p$ and each $\gamma_{i}$ is positive in $p$.
        \end{enumerate}
\end{enumerate}

\paragraph{Phase 3: Translation and Output.}  If there is some system in the Reduction and Elimination stage from some propositional variable could be eliminated, then ALBA reports failure and terminates. Otherwise, every initial system $\{\i_{0} \leq \phi_{i}~,~\i_{0}\leq\mathbf{a}~,~ \psi_{i} \leq \m_{0}\}$ has been reduced to a system, denoted $\mathit{Reduce}(\phi_{i} \wedge \mathbf{a}\leq \psi_{i})$, containing no propositional variables. Let $ALBA(\phi \wedge \mathbf{a}\leq \psi)$ be the set of quasi-inequalities: 
\[
\& \left( \mathit{Reduce}(\phi_{i} \wedge \mathbf{a}\leq \psi_{i}) \right)\Rightarrow \i_{0} \leq \m_{0}
\]
for each $\phi_{i} \wedge \mathbf{a}\leq \psi_{i} \in \Pre(\phi \wedge \mathbf{a}\leq \psi)$.  All members of $ALBA(\phi \wedge \mathbf{a}\leq \psi)$ are free of propositional variables, so applying the standard translation to each member of $ALBA(\phi \wedge \mathbf{a}\leq \psi)$ will result in a set of first-order correspondents -- one for each member of the set of quasi-inequalities.  Let $ALBA^{FO}(\phi \wedge \mathbf{a}\leq \psi)$ equal:
%
%\begin{eqnarray*}
%    \bigwedge_{1\leq i\leq n}\forall \bC_{\i_{0}} \forall \bc_{\m_{0}}\forall \bC_{\m_{0}} \forall\overline{\mathbf{c}_{\j}}\forall\overline{\mathbf{C}_{\j}}\forall \overline{\bc_{\m}}\forall \overline{\bC_{\m}}  \Big( \forall x ST_{x} \Big(\bigwedge {{Reduce}}(\phi_{i} \wedge \mathbf{a} \leq \psi_{i})\Big) &\Rightarrow& \forall x ST_{x}(\i_{0} \leq \m_{0})\Big)
%\end{eqnarray*}
%
%\begin{eqnarray*}
%    \bigwedge_{1\leq i\leq n}\forall \bC_{\i_{0}} \forall \bc_{\m_{0}}\forall \bC_{\m_{0}} \forall\overline{\mathbf{c}_{\j}}\forall\overline{\mathbf{C}_{\j}}\forall \overline{\bc_{\m}}\forall \overline{\bC_{\m}} \Big( \bigwedge \{ \ST_x(\alpha \leq \beta) \mid \alpha \leq \beta \in {{Reduce}}(\phi_{i} \wedge \mathbf{a} \leq \psi_{i})\} \to \forall x ST_{x}(\i_{0} \leq \m_{0})\Big)
%\end{eqnarray*}

\begin{eqnarray*}
    \bigwedge_{1\leq i\leq n}\forall \bC_{\i_{0}} \forall \bc_{\m_{0}}\forall \bC_{\m_{0}} \forall\overline{\mathbf{c}_{\j}}\forall\overline{\mathbf{C}_{\j}}\forall \overline{\bc_{\m}}\forall \overline{\bC_{\m}} \Big( \bigwedge \ST_x({{Reduce}}(\phi_{i} \wedge \mathbf{a} \leq \psi_{i})) \to \forall x ST_{x}(\i_{0} \leq \m_{0})\Big)
\end{eqnarray*}
where we write $\ST_x({{Reduce}}(\phi_{i} \wedge \mathbf{a} \leq \psi_{i})$ for $\{ \ST_x(\alpha \leq \beta) \mid \alpha \leq \beta \in {{Reduce}}(\phi_{i} \wedge \mathbf{a} \leq \psi_{i})\}$ and $\overline{\mathbf{c}_{\j}}$, $\overline{\mathbf{C}_{\j}}$, $\overline{\mathbf{c}_{\m}}$  and $\overline{\mathbf{C}_{\m}}$ are the vectors of all variables corresponding to the standard translations of nominals and conominals, respectively, occurring in $\mathit{Reduce}(\phi_{i} \wedge \mathbf{a}\leq \psi_{i})$, other than the reserved nominal and conominal $\i_{0}$ and $\m_{0}$.  Note that the constant $\mathbf{c}_{\i_{0}}$ corresponding to $\i_{0}$ is not quantified.  This is done to produce a \emph{local} frame correspondent. For a global correspondent, simply take the universal closure. In practice, as illustrated in Example \ref{examp:alba:reflex} below, it often helps to simplify the quasi-inequalities ${{Reduce}}(\phi_{i} \wedge \mathbf{a} \leq \psi_{i})$ before applying the standard translation. 

\begin{example}\label{examp:alba:reflex}
Let $\mathbb{A}$ be an arbitrary finite Heyting algebra and let $a\in A$.  As already discussed in Section \ref{sec:mv:corresp:theory}, it is established in \cite{Frank} that $p\rightarrow \Diamond p$ and $a$-reflexivity are local frame $a$-correspondents.  We now show how ALBA can be used to obtain the same result. 

The inequality $p \wedge \mathbf{a} \leq \Diamond p$ requires no preprocessing, hence the first-approximation rule may be applied to it to obtain $\{\nomi_{0} \leq p \wedge \mathbf{a} ~,~ \Diamond p\leq \cnomm_{0}\}$.  The splitting rule can immediately be applied to this system to obtain $\{\nomi_{0} \leq p~,~\nomi_{0} \leq \mathbf{a} ~,~ \Diamond p\leq \cnomm_{0}\}$.  Applying the right-hand Ackermann rule, the system $\{\nomi_{0} \leq \mathbf{a} ~,~ \Diamond\nomi_{0}\leq \cnomm_{0}\}$ is obtained.  $ALBA^{FO}(p \wedge \mathbf{a} \leq \Diamond p)$ then equals:
\begin{eqnarray*}
 &\quad& \forall \mathbf{C}_{\nomi_{0}}\forall \mathbf{C}_{\cnomm_{0}}\forall \mathbf{c}_{\cnomm_{0}}\Big[\forall x\Big[\Big((x= \mathbf{c}_{\nomi_{0}}\wedge \mathbf{C}_{\nomi_{0}}) \preceq \mathbf{a}\Big) 
\wedge \Big(\exists y\Big(Rxy \wedge (y= \mathbf{c}_{\nomi_{0}}\wedge \mathbf{C}_{\nomi_{0}})\Big) \preceq (x\neq \mathbf{c}_{\cnomm_{0}}\vee \mathbf{C}_{\cnomm_{0}})\Big)\Big]\\
&\Rightarrow& \forall x\Big((x= \mathbf{c}_{\nomi_{0}}\wedge \mathbf{C}_{\nomi_{0}}) \preceq (x\neq \mathbf{c}_{\cnomm_{0}}\vee \mathbf{C}_{\cnomm_{0}})\Big)\Big].
\label{ALBAreflexivity2}
\end{eqnarray*} 
To see how this may be simplified, consider an assignment $v$ on an $\mathbb{A}$-frame considered as a first-order structure. Under $v$, $ALBA^{FO}(p \wedge \mathbf{a} \leq \Diamond p)$ is equivalent to:
\begin{eqnarray*}
&\quad&\bigwedge_{j\in J^{\infty}(\mathbb{A})}\bigwedge_{m\in M^{\infty}(\mathbb{A})}\bigwedge_{c_{\nomi_{0}}\in W}\Big[\bigwedge_{x\in W}\Big[\Big(\Big((v(x)= c_{\nomi_{0}})\wedge j\Big) \leq a\Big) 
\wedge \Big(\bigvee_{y\in W}\Big(v(R(x,y)) \wedge \Big((v(y)= c_{\nomi_{0}})\wedge j\Big)\Big)\\
&\leq& \Big((v(x)\neq c_{\cnomm_{0}})\vee m\Big)\Big)\Big]\Rightarrow\bigwedge_{x\in W}\Big(\Big((v(x)= c_{\nomi_{0}})\wedge j\Big) \leq \Big((v(x)\neq c_{\cnomm_{0}})\vee m\Big)\Big)\Big]
\label{eq:}
\end{eqnarray*}
which simplifies to:
\[
\bigwedge_{j\in J^{\infty}(\mathbb{A})}\bigwedge_{m\in M^{\infty}(\mathbb{A})}\Big(\Big((j \leq a) \wedge ((R(c_{\nomi_{0}},c_{\nomi_{0}}) \wedge j) \leq m)\Big) \Rightarrow  (j \leq m)\Big),
\]
which, in turn, is equivalent to $(\forall m \in M^{\infty}(\mathbb{A}))\Big( a \wedge  (R(c_{\nomi_{0}},c_{\nomi_{0}}) \to m) \leq m \Big)$, which may be further simplified as follows:
%
%\begin{eqnarray}
%    &&(\forall m \in M^{\infty}(\mathbb{A}))\Big( a \wedge  (R(c_{\nomi_{0}},c_{\nomi_{0}}) \to m) \leq m \Big)\nonumber\\
    %
%    &\text{iff} &(\forall m \in M^{\infty}(\mathbb{A}))\Big( \text{if } R(c_{\nomi_{0}},c_{\nomi_{0}}) \to m) \not\leq m  \:\text{ then }\: a \leq m \Big)\label{eq:m:prime:m}\\
    %
%    &\text{iff} &(\forall m \in M^{\infty}(\mathbb{A}))\Big( \text{if } \lambda(m) \leq (R(c_{\nomi_{0}},c_{\nomi_{0}}) \to m) \:\text{ then }\: a \leq m \Big)\label{eq:apply:lambda}\\
    %
%    &\text{iff} &(\forall m \in M^{\infty}(\mathbb{A}))\Big( \text{if } \lambda(m) \wedge R(c_{\nomi_{0}},c_{\nomi_{0}}) \leq  m \:\text{ then }\: a \leq m \Big)\label{eq:apply:residuation}\\
    %
%    &\text{iff} &(\forall m \in M^{\infty}(\mathbb{A}))\Big( \text{if } R(c_{\nomi_{0}},c_{\nomi_{0}}) \leq  m \:\text{ then }\: a \leq m \Big)\label{eq:m:not:above:lambda:m}\\
    %
%    &\text{iff} &a \leq R(c_{\nomi_{0}},c_{\nomi_{0}})\label{eq:meet:generation}
%\end{eqnarray}
%
%
\begin{eqnarray}
    &&a \wedge  (R(c_{\nomi_{0}},c_{\nomi_{0}}) \to m) \leq m,  \quad \quad \text{for all } m \in M^{\infty}(\mathbb{A}))\nonumber\\
    &\text{iff} &\text{if } R(c_{\nomi_{0}},c_{\nomi_{0}}) \to m) \not\leq m  \:\text{ then }\: a \leq m,  \quad \quad \text{for all } m \in M^{\infty}(\mathbb{A})) \label{eq:m:prime:m}\\
    &\text{iff} &\text{if } \lambda(m) \leq (R(c_{\nomi_{0}},c_{\nomi_{0}}) \to m) \:\text{ then }\: a \leq m,  \quad \quad \text{for all } m \in M^{\infty}(\mathbb{A}))\label{eq:apply:lambda}\\
    &\text{iff} &\text{if } \lambda(m) \wedge R(c_{\nomi_{0}},c_{\nomi_{0}}) \leq  m \:\text{ then }\: a \leq m,  \quad \quad \text{for all } m \in M^{\infty}(\mathbb{A}))\label{eq:apply:residuation}\\
    &\text{iff} &\text{if } R(c_{\nomi_{0}},c_{\nomi_{0}}) \leq  m \:\text{ then }\: a \leq m,  \quad \quad \text{for all } m \in M^{\infty}(\mathbb{A}))\label{eq:m:not:above:lambda:m}\\
    &\text{iff} &a \leq R(c_{\nomi_{0}},c_{\nomi_{0}})\label{eq:meet:generation}
\end{eqnarray}
where \eqref{eq:m:prime:m} follows by the meet-primeness  of $m$, \eqref{eq:apply:lambda} follows by \eqref{eq:lambda}, \eqref{eq:apply:residuation} by residuation, \eqref{eq:m:not:above:lambda:m} by the meet-primeness of $m$ and the fact that $\lambda(m) \not\leq m$, and lastly \eqref{eq:meet:generation} follows by the complete meet-generation of $\algA$ by $\mty(\algA)$.
\end{example}

\subsection{Correctness of ALBA and completeness for inductive formulas}

In this subsection, we give a sketch of the correctness of the algorithm described above, thereby establishing that, when ALBA succeeds on an $\mvml$-formula or inequality, the first-order formula produced is indeed its $a$-frame correspondent. We also argue that ALBA succeeds on all inductive $\mvml$-formulas, and hence that they all have effectively computable $a$-frame correspondents.

%Note that the condition $v(\phi)\wedge \mathbf{a} \leq v(\psi)$ is equivalent to $v(\mathbf{a}) \leq v(\phi)\rightarrow v(\psi)$ under all assignments $v$ on $\mathbb{G}$.

\begin{thm}
If ALBA succeeds in reducing an inequality $\phi \leq \psi$, then, for any $\mathbb{A}$-frame $\modF = (W,R)$, 
\begin{eqnarray*}
&\quad& \mathfrak{F},w\Vdash_{a} \phi\leq\psi\\
&\mbox{ iff }& \mathfrak{F}^{+} (\i_{0}:= f_{\i_{0}})\models ALBA(\phi \wedge \mathbf{a}\leq \psi) \\
&\mbox{ iff }& \mathfrak{F} \models ALBA^{FO}(\phi \wedge \mathbf{a}\leq \psi)[c_{\i_{0}}:=w]
\label{eq:}
\end{eqnarray*}
where $f_{\i_{0}}\in A^{W}$ such that $f_{\i_{0}}(w)\in J^{\infty}(\mathbb{A})$ for exactly one $w\in W$ and $f_{\i_{0}}(w_{0})=0$ for all $w_{0}\neq w$.
\label{ALBACorrect2}
\end{thm}

\begin{proof}
It is sufficient to show the following:
\begin{eqnarray}
&\quad& \modF , w \Vdash_{a} \phi \leq \psi\label{eq:correctness:1}\\
&\mbox{ iff }&\mathfrak{F},w \Vdash  \phi \wedge \ba\leq \psi \label{eq:correctness:2}\\
%&\mbox{ iff }& \mathfrak{F}^{+} (\i_{0}:= f_{\i_{0}})\models_{a} \phi \leq \psi \\
%&\mbox{ iff }& \mathfrak{F}^{+} (\i_{0}:= f_{\i_{0}})\models \phi\wedge \mathbf{a} \leq \psi \\
&\mbox{ iff }& \modF, w \Vdash \Pre(\phi \wedge \ba \leq \psi)\label{eq:correctness:3}\\
&\mbox{ iff }&\modF, V \Vdash (\i_0 \leq \phi_{i} \wedge \ba~\&~ \psi_{i}\leq \m_{0}) \Rightarrow (\i_{0} \leq  \m_{0}), \quad V(\nomi_0, w) \in \jty(\algA), \: 1 \leq i \leq \ell\label{eq:correctness:4}\\
&\mbox{ iff }&\modF, V \Vdash \mathit{Reduce}(\phi_{i} \wedge \ba \leq  \psi_{i}) \Rightarrow (\i_{0} \leq \m_{0}), \quad V(\nomi_0, w) \in \jty(\algA), \: 1 \leq i \leq \ell\label{eq:correctness:5}\\
&\mbox{ iff }&\modF, V \Vdash ALBA(\phi \wedge \mathbf{a}\leq \psi), \quad V(\nomi_0, w) \in \jty(\algA)\label{eq:correctness:6}\\
&\mbox{ iff }&\mathfrak{F} \models ALBA^{FO}(\phi \wedge\mathbf{a}\leq \psi)[c_{\i_{0}}:=w]
\label{eq:correctness:7}
\end{eqnarray}
\noindent
\underline{\eqref{eq:correctness:1} iff \eqref{eq:correctness:2}:} By definition of $a$-truth of inequalities. \\

\noindent
\underline{\eqref{eq:correctness:2} iff \eqref{eq:correctness:3}:} It is easy to check that the preprocessing steps preserve validity.  The rules for the elimination of monotone and antitone variables follow by the natural many-valued analogue of the monotonicity (antitonicity) of positive (negative) formulas. That is, if $\phi$ is positive (negative) in $p$, then for any valuations $V$ and $V'$ on a frame $\modF = (W,R)$ such that $V(p,w) \leq V'(p,w)$ for all $w \in W$, it holds that $V(\phi,w) \leq V'(\phi,w)$ ($V'(\phi,w) \leq V(\phi,w)$) for all $w \in W$.\\

\noindent
\underline{\eqref{eq:correctness:3} iff \eqref{eq:correctness:4}:} This equivalence requires some justification. Particularly, we have to justify the local validity expressed in \eqref{eq:correctness:3} being transformed into an expression of global validity in \eqref{eq:correctness:4}. 

Assume that $\modF, w \Vdash \Pre(\phi \wedge \ba \leq \psi)$, i.e., $\modF, w \Vdash \phi_i \wedge \ba \leq \psi_i$ for each  $\phi_{i} \wedge \mathbf{a}\leq \psi_{i} \in \Pre(\phi  \wedge \ba\leq \psi)$. Suppose that $V$ is a valuation on $\modF$ such that $V(\nomi_0, w) = i_0 \in \jty(\algA)$ and $\modF, V \Vdash \nomi_0 \leq \phi_i \wedge \ba$ and $\modF, V \Vdash \psi_i \leq \cnomm_0$. But then, by assumption, $V(\nomi_0, w) \leq V(\phi_i \wedge \ba, w) \leq V(\psi_i, w) \leq V(\cnomm_0, w)$. Moreover, for states $v \neq w$, we have $V(\nomi_0, v) = 0 \leq V(\cnomm_0, v)$. Therefore, $\modF, V \Vdash \nomi_0 \leq \cnomm_0$.

Conversely, assume that $\modF, V \Vdash (\i_0 \leq \phi_{i} \wedge \ba~\&~ \psi_{i}\leq \m_{0}) \Rightarrow (\i_{0} \leq  \m_{0})$ for all valuations $V$ such that  $V(\nomi_0, w) \in \jty(\algA)$ and all $1 \leq i \leq \ell$. Let $V$ be any valuation on $\modF$. Then $V(\phi_i \wedge \ba, w) = \bigvee \{i \in \jty(\algA) \mid i \leq V(\phi_i \wedge \ba, w) \}$ and $V(\psi_i, w) = \bigwedge \{m \in \mty(\algA) \mid V(\psi_i, w) \leq m \}$. Thus to show that $V(\phi_i \wedge \ba, w) \leq V(\psi_i, w)$, it is sufficient to show that for every $i \in \jty(\algA)$ with  $i \leq V(\phi_i \wedge \ba, w)$ and every $m \in \mty(\algA)$ with  $V(\psi_i, w) \leq m$, it is the case that $i \leq m$.  Now, if $i \leq V(\phi_i \wedge \ba, w)$ and $V(\psi_i, w) \leq m$, let $V'$ be the valuation which is just like $V$ except that $V'(\nomi_{0}, w) = i$ and $V(\cnomm_0, w) = m$. Since $\nomi_0$ and $\cnomm_0$ occur in neither $\phi_i \wedge \ba$ nor $\psi_i$, it is the case that $V'(\nomi_{0}, w) \leq V'(\phi_i \wedge \ba, w)$ and $V'(\psi_i, w) \leq V(\cnomm_0, w)$ and so, by \eqref{eq:correctness:4}, $i = V'(\nomi_{0}, w) \leq V'(\cnomm_0, w) = m$, as desired. \\

\noindent
\underline{\eqref{eq:correctness:4} iff \eqref{eq:correctness:5}:} This follows by the soundness, as proven in \cite{ALBA}, of the rules of the reduction and elimination phase on perfect distribute lattice expansion, including perfect Heyting algebras like $\modF^+$.\\

\noindent
\underline{\eqref{eq:correctness:5} iff \eqref{eq:correctness:6}:} A trivial consequence of the definition of  $ALBA(\phi \wedge \mathbf{a}\leq \psi)$, as given in the description of phase 3, above. 

\noindent
\underline{\eqref{eq:correctness:6} iff \eqref{eq:correctness:7}:} This follows from Corollary \ref{Corollary:StandardTranslation}.
\end{proof}

In \cite[Theorem 8.1]{ALBA} and \cite[Theorem 8.8]{nonDistALBA} it is proved that ALBA succeeds on all inductive inequalities for logics with algebraic semantics given by distributive, respectively general, lattice expansions. We have indicated above that the only difference between an execution of ALBA for $a$-validity and an execution for validity (as in \cite{ALBA}), is that the inequality $\i_0 \leq \mathbf{a}$ is being carried along. Since this inequality is effectively ``inert'' during the entire execution, the proof of theorem Theorem 8.1 in \cite{ALBA} can be applied essentially unchanged. We therefore have the following theorem:

\begin{thm}
ALBA succeeds on the class of all inductive inequalities (and hence the class of all Sahlqvist inequalities) of $L^{+}_{\mathbb{A}}$ over all $\mathbb{A}$-frames.
\label{ALBATheorem}
\end{thm}

%Moreover, the proof of correctness of ALBA given in \cite{ALBA} is a special case of Theorem \ref{ALBATheorem} with $a=1$.  In other words, suppose that $\phi\leq \psi$ is an inductive inequality or a Sahlqvist inequality.  Then ALBA succeeds on $\phi\leq\psi$ and produces $ALBA^{FO}(\phi\leq\psi)$.  If we start with the inequality $\phi\wedge \mathbf{a}\leq \psi$, with $a=1$, then $\phi\wedge \mathbf{1}\leq \psi$ will be reduced to $ALBA^{FO}(\phi\wedge\mathbf{1}\leq\psi)=ALBA^{FO}(\phi\leq\psi)$ by ALBA.
%
Combining Theorems \ref{ALBACorrect2} and \ref{ALBATheorem} yields the following corollary:

\begin{cor}
All inductive formulas have effectively computable first-order local frame $a$-correspondents over all $\mathbb{A}$-frames.
\label{MajorCorR}
\end{cor}
%We now have a sufficient condition that assures us that an inequality will have a first-order frame correspondent, and that condition is that the inequality be inductive.  This is a very strong Sahlqvist type result, as it is applicable to a wide range of truth-value spaces.\wmnote{Rephrase.}

\section{Comparing classical and many-valued frame correspondents}\label{Section:SvBResult}
In Section \ref{Section:ALBA} we saw that each inductive formula, and hence each Sahlqvist formula, has a first-order $a$-correspondent on frames. Moreover, it already follows from the work of Frankowski \cite{Frank}, as discussed in Section \ref{sec:mv:corresp:theory}, that for some of the most well-known modal axioms (including reflexivity, transitivity, symmetry, density and right-unboundedness) the first-order $a$-correspondents are syntactically identical to the classical ones. Indeed, as we will now see, this is an instance of a general phenomenon: in this section we show that the $a$-correspondents for Sahlqvist formulas are \emph{syntactically identical} to the classical correspondents produced by the well-known Sahlqvist-van Benthem algorithm.  We do this by introducing an analogous algorithm for the present, many-valued setting. Although this produces the same output as the classical algorithm, the inner workings are complicated by the fact that classical minimal valuations are no longer minimal.

Let $\sigma: \mvsol \to \emvsol$ be a predicate substitution acting on predicate symbol $\P$ given by $\lambda u. \delta(u)$ for some $\emvsol$-formula $\delta(u)$, i.e. when applied to a formula $\phi \in \mvsol$, $\sigma(\phi)$ is the formula obtained by replacing each occurrence of an atomic subformula $\P(t)$ with $\delta(t)$. If $\bC$ is a truth-value constant/variable, the substitution $\sigma_{\bC}$ is like $\sigma$ but replaces each occurrence of an atomic subformula $\P(t)$ with $\delta(t) \wedge \bC$. 

%\marginpar{Introduce the notation and talk about the substitution instances.}
%\begin{align*}
%\sigma_\bC(P(y))&=\left(\delta(y)\wedge \bC\right)\\
%\sigma(P(y))&=\delta(y)
%\end{align*}
%where $\delta(y)$ is some function of $y$ and $C\in J^{\infty}(\algA)$ corresponds to $\bC$.

\begin{lemma}\label{Lemma:EliminatingC}
	Let $\phi$ be any formula of $\mvsol$.	Then for any $\algA$-model, $\modF$, for $\mvsol$ and any assignment $v$ on $\modF$, it holds that
	\[
	\Vert\sigma_{\bC}(\phi) \wedge \bC\Vert_{\modF,v} = \Vert\sigma(\phi) \wedge \bC\Vert_{\modF,v}\\
	\]
\end{lemma}

\begin{proof}
	The proof proceeds by induction on $\phi$. The base cases for atomic formulas which are truth value constants or truth value variables are straightforward. Consider the case when $\phi$ is an atomic formula of the form $\P(t)$, then 
	\[
	\Vert\sigma_{\bC}(\phi)\wedge \bC\Vert_{\modF,v}=
	\Vert\sigma_{\bC}(\P(t))\wedge \bC\Vert_{\modF,v}=
	\Vert\left(\delta(t)\wedge \bC\right)\wedge \bC\Vert_{\modF,v}=
	\Vert \delta(t) \wedge \bC\Vert_{\modF,v}=
	\Vert\sigma(\phi)\wedge \bC\Vert_{\modF,v}.
	\]
	\noindent If $\phi$ is an atomic formula of the form $\P'(t)$ where $\P'$ is different from $\P$, then 
	\[
	\Vert\sigma_{\bC}(\phi)\wedge \bC\Vert_{\modF,v}=
	\Vert\sigma_{\bC}(\P'(t))\wedge \bC\Vert_{\modF,v}=
	\Vert \P'(t) \wedge \bC \Vert_{\modF,v}=
	\Vert \sigma(\P'(t)) \wedge \bC\Vert_{\modF,v}=
	\Vert\sigma(\phi)\wedge \bC\Vert_{\modF,v}.
	\]
	\noindent The cases when $\phi$ is an atomic formula of the form $t_1 = t_2$ or $Rt_1t_2$ are similar. 
	%If $\phi$ is an atomic formula of the form $t_1 = t_2$, then 
	%\begin{align*}
	%	\Vert\sigma_{\bC}(\phi)\wedge \bC\Vert_{\modF,v} &= \Vert\sigma_{\bC}(t_1 = t_2)\wedge \bC\Vert_{\modF,v}\\
	%		&= \Vert ((t_1 = t_2) \wedge \bC) \wedge \bC\Vert_{\modF,v}\\
	%	&= \Vert (t_1 = t_2) \wedge \bC\Vert_{\modF,v}\\
	%	&= \Vert\sigma(\phi)\wedge \bC\Vert_{\modF,v}.
	%\end{align*}
	%	
	If $\phi$ is of the form $\psi \star \chi$ for $\star\in\{\vee,\wedge\}$, then using distributivity if $\star$ is $\vee$ and associativity and indempotency if $\star$ is $\wedge$ we have that:
	\begin{align*}
		\Vert\sigma_{\bC}(\psi \star \chi)\wedge \bC\Vert_{\modF,v}
		&=\Vert(\sigma_{\bC}(\psi) \star \sigma_{\bC}(\chi)) \wedge \bC\Vert_{\modF,v}\\
		&=\Vert \sigma_{\bC}(\psi) \wedge \bC \Vert_{\modF,v} \star^{\mathbb A} \Vert \sigma_{\bC}(\chi) \wedge \bC \Vert_{\modF,v}\\
		%&=\Vert\sigma_{\bC}(\psi \wedge \bC) \star \sigma_{\bC}(\chi)\Vert_{\modF,v}\\
		%&=\left(\Vert\sigma_{\bC}(\psi)\Vert_{\modF,v} \star \Vert\sigma_{\bC}(\chi)\Vert_{\modF,v}\right)\wedge C\\
		&= \Vert\sigma(\psi)\wedge \bC\Vert_{\modF,v} 
		\star^{\mathbb A} \Vert\sigma(\chi)\wedge \bC\Vert_{\modF,v} &&\text{(by the I.H.)}\\
		&= \Vert (\sigma(\psi)\wedge \bC) \star (\sigma(\chi)\wedge \bC)\Vert_{\modF,v}\\
		&= \Vert\sigma(\psi \star \chi)\wedge \bC\Vert_{\modF,v}\\
	\end{align*}
	
	\noindent Next suppose that $\phi$ is of the form $\psi\rightarrow\chi$.  We will use the Heyting algebra identity 
	$(a\rightarrow b)\wedge c=((a\wedge c)\rightarrow (b\wedge c))\wedge c$.  Here
	\begin{align*}
		\Vert\sigma_{\mathbf{C}}(\psi \rightarrow \chi)\wedge \bC\Vert_{\modF,v}
		&=\Vert(\sigma_{\mathbf{C}}(\psi) \rightarrow \sigma_{\mathbf{C}}(\chi))\wedge \bC\Vert_{\modF,v}\\
		&=\Vert((\sigma_{\mathbf{C}}(\psi) \wedge \bC) \rightarrow (\sigma_{\mathbf{C}}(\chi) \wedge \bC))\wedge \bC\Vert_{\modF,v} &&(\text{Heyting alg.\ indentity})\\
		&=\left(\Vert\sigma_{\bC}(\psi) \wedge \bC \Vert_{\modF,v} \rightarrow^{\mathbb A} \Vert\sigma_{\bC}(\chi \wedge \bC)\Vert_{\modF,v}\right) \wedge^{\mathbb A} \Vert\bC\Vert_{\modF,v}\\
		&= \left(\Vert\sigma(\psi)\wedge \bC\Vert_{\modF,v} \rightarrow^{\mathbb A} 
		\Vert\sigma(\chi)\wedge \bC\Vert_{\modF,v}\right)\wedge^{\mathbb A} \Vert\bC\Vert_{\modF,v} &&\text{(by the I.H.)}\\
		&= \Vert \left(\sigma(\psi)\wedge \bC \rightarrow \sigma(\chi) \wedge \bC \right) \wedge \bC\Vert_{\modF,v}\\
		&= \Vert\sigma(\psi\rightarrow\chi)\wedge\bC\Vert_{\modF,v} &&(\text{Heyting alg.\ indentity}).
	\end{align*}
	
	\noindent Suppose $\phi$ is of the form $\exists x \psi$.  Then
	\begin{align*}
		\Vert\sigma_{\bC}(\exists x \psi)\wedge \bC\Vert_{\modF,v}
		&=\bigvee \left\{\Vert\sigma_{\bC}(\psi)\Vert_{\modF,v'}\mid v'\sim_{x}v\right\}\wedge^{\mathbb A} \Vert \bC \Vert_{\modF,v}\\
		&= \bigvee \left\{\Vert\sigma_{\bC}(\psi)\wedge \bC\Vert_{\modF,v'}\mid v'\sim_{x}v\right\}\\
		&= \bigvee \left\{\Vert\sigma(\psi)\wedge \bC\Vert_{\modF,v'}\mid v'\sim_{x}v\right\} &&\text{(by the I.H.)}\\
		&=\bigvee \left\{\Vert\sigma(\psi)\Vert_{\modF,v'}\mid v'\sim_{x}v\right\}\wedge^{\mathbb A} \Vert \bC \Vert_{\modF,v}\\
		%&=\bigvee \left\{\Vert\sigma(\psi)\Vert_{\modF,v}\mid v'\sim_{x}v\right\} \wedge \Vert\bC\Vert_{\modF,v}\wedge C\\
		%&=\Vert\sigma(\exists x \psi)\Vert_{\modF,v} \wedge \Vert\bC\Vert_{\modF,v}\\
		&=\Vert\sigma(\exists x \psi) \wedge \bC\Vert_{\modF,v}.
	\end{align*}
	
	\noindent The case for $\phi$ of the form $\forall x \psi$ is similar.
\end{proof}

\noindent The following lemma collects the $\mathbb{A}$-valued analogues of some standard first-order equivalences which will be used below. 

\begin{lemma}\label{Lemma:Equivalences}
If $\alpha$, $\beta$ and $\gamma$ are formulas of $\emvfol$ or $\emvsol$, then:
\begin{enumerate}[label=(\roman*)]
\item $\left(\exists x \alpha(x) \wedge \beta\right)\equiv \exists x \left(\alpha(x) \wedge \beta \right)$
if the variable $x$ does not occur free in $\beta$.
\item $\left(\exists x \alpha(x) \rightarrow \beta\right)\equiv \forall x \left(\alpha(x) \rightarrow \beta \right)$
if the variable $x$ does not occur free in $\beta$.
\item $\left((\alpha \vee \beta)\rightarrow \gamma\right)\equiv 
\left((\alpha \rightarrow \gamma)\wedge(\beta \rightarrow \gamma)\right).$
\item $\forall x_{1},x_{2},\hdots,x_{n} \left(\alpha \wedge \beta\right)\equiv 
\left(\forall x_{1},x_{2},\hdots,x_{n} \alpha \wedge \forall x_{1},x_{2},\hdots,x_{n} \beta \right).$
\end{enumerate}
\end{lemma}

We will now recall the classical definition of Sahlqvist formulas. We will refer to these formulas as the {\em classical Sahlqvist formulas}, to distinguish them from the formulas of Definition \ref{Sahlqvist:Ineq:Def}. The relationship between these two classes will be clarified below. 

\begin{definition}[Classical Sahlqvist formulas (see e.g.\ \cite{vanb, BlueBook})]\label{Def:Classical:Sahl}
	A {\em boxed atom} is a formula $\Box^n p$ where $\Box^n$ is a string of $n$ copies of $\Box$ for $n \geq 0$. A
	{\em (definite) classical Sahlqvist antecedent} is a formula built up from $\bo$, $\mathbf{1}$, boxed
	atoms and negative formulas using $\wedge,\vee$ and $\Diamond$ ($\wedge$ and $\Diamond$).  A {\em (definite) classical Sahlqvist implication}
	is an implication $\phi\rightarrow \psi$ such that $\psi$ is positive and $\phi$ a (definite) classical Sahlqvist 
	antecedent. A {\em (definite) classical Sahlqvist formula} is a formula that is built up from (definite) classical Sahlqvist implications using
	$\Box$ and $\wedge$ without restriction, and using $\vee$ between formulas with no proposition variables in common.
\end{definition}

The following Lemma is an easy exercise (see e.g.\ \cite{conradie2017algebraic}):

\begin{lemma}\label{lem:Definite:Classical:Sahl}
	Every classical Sahlqvist implication is equivalent to a conjunction of definite classical Sahlqvist implications.
\end{lemma}

It is not difficult to see that if $\phi$ is a classical Sahlqvist formula and we read classical negations $\neg \psi$ as $\psi \to \bot$, then $\top \leq \phi$ is $\epsilon$-Sahlqvist for $\epsilon$ constantly $1$. Thus, when projecting onto the classical modal logic case, Definition \ref{Sahlqvist:Form:Def} covers all classical Sahlqvist formulas. On the purely syntactic level, this projection goes beyond Defintion \ref{Def:Classical:Sahl}. For example, $\top \leq \Box \neg p \to \Box \Box \neg p$ is Sahlqvist but $\Box \neg p \to \Box \Box \neg p$ is not a classical Sahlqvist formula. The latter is, however, locally equivalent on frames to the well-known Sahlqvist formula $\Box p \to \Box \Box p$, obtained by switching the polarity of the proposition letter $p$. The class of {\em Sahlqvist-van Benthem formulas} \cite{kracht1999tools} enlarges the class of classical Sahlqvist formulas to accommodate such polarity switching. For a more systematic comparison we refer the reader to \cite[Section 3.3]{ALBA} which also provides, although in a slightly different context, the arguments needed for the proof of the following proposition:

\begin{prop}
	Every Sahlqvist formula is locally equivalent over $\mathbf{2}$-frames to a classical Sahlqvist formula and vice versa.
\end{prop}

Lemma \ref{lem:Definite:Classical:Sahl} can be used in conjunction with the following two lemmas to reduce finding local $a$-correspondents for classical Sahlqvist formulas to finding local $a$-correspondents for classical Sahlqvist implications.

\begin{lemma}\label{Lemma:Disjunctions}%\wmnote{Possibly this lemma needs to be moved as well.}
Let $\modF = (W,R)$ be an $\algA$-frame, $w \in W$ and $\phi, \psi \in \mvml$-formulas that
do not share variables.  Then, $\modF, w \Vdash_{a} \phi \vee \psi$ if, and only if, there are $a_1, a_2 \in A$
such that  $\modF, w \Vdash_{a_1} \phi$ and  $\modF, w \Vdash_{a_2} \psi$ and $a_1 \vee a_2 \geq a$.
\end{lemma}

\begin{proof}
Suppose that $\modF, w \Vdash_{a} \phi \vee \psi$. Let $S = \{ V(\phi, w) \mid V \text{ a valuation on } \modF \}$ and 
$T = \{ V(\psi, w) \mid V \text{ a valuation on } \modF \}$. Then, by
the assumptions that $\phi$ and $\psi$ share no variables and that
$\modF, w \Vdash_{a} \phi \vee \psi$, it follows that $s \vee t \geq a$
for all $s \in S$ and $t \in T$. Therefore
$a \leq \bigwedge \{ s \vee t \mid s \in S \text{ and } t \in T\}  =
\bigwedge T \vee \bigwedge S$. We can therefore take $b = \bigwedge T$ and $c = \bigwedge S$. The other direction is immediate.
\end{proof}

\noindent For the sake of translating boxed atoms into $\mvfol$ succinctly we will, as normal, write $R^{n}xy$ for $\exists z_1 \cdots \exists z_n(Rxz_1 \wedge Rz_1z_2 \wedge \cdots \wedge Rz_{n-1}z_n \wedge z_n = y)$ when $n \geq 1$, and $R^{0}xy$ for $x = y$. We will write $\alpha[y/x]$ for the $\mvfol$ (or $\mvsol$) formula obtained from $\alpha$ by renaming all free occurrences of $x$ to $y$ or simply $\alpha(y)$ if $x$ is clear from the context. A straightforward proof by induction on $n$ suffices to establish the following lemma.

\begin{lemma}\label{lem:translate:boxed:atom}
    For any $\emvml$-formula of the form $\Box^{n} \phi$, it is the case that $\ST_{x}(\Box^{n} \phi) \equiv \forall y (R^n xy \to \ST_{y}(\phi))$. 
\end{lemma}

%\begin{proof}
%One may prove this lemma rigorously proceeding by induction on $n$. However, teh gist is already  but contains the key steps  FOr the sake of berWe do the proof for the case $n$
%If $n = 0$ then $\ST_{x}(\Box^{n} \phi) \equiv \ST_{x}(\phi) \equiv \forall y (x = y \to \ST_{y}(\phi)) \equiv \forall y (R^n xy \to \ST_{y}(\phi))$. Assuming the claim holds for $n = k$, consider  $\ST_{x}(\Box^{k+1} \phi)$. Now $\ST_{x}(\Box^{k+1} \phi) \equiv \ST_{x}(\Box^{k}\Box \phi) \equiv \forall u (R^kxu \to \ST_{u}(\Box \phi) \equiv $. 
%\end{proof}

\begin{lemma}\label{Lemma:SahlFormulasToImplications}%\wmnote{The clause for $\vee$ in this lemma did not work as it was previously formulated. Also, the clause for $\Box$ needed stronger assumptions and a new proof.}
	Let $\phi$ and $\psi$ be formulas of $\mvml$ and $a\in A$.
	\begin{enumerate}[label=(\roman*)]
		\item If $\phi$ and $\alpha(x)$ are local frame $b$-correspondents for all $b \in A$,
		then $\Box^{n} \phi$ and $\forall y \left(R^{n}xy \rightarrow \alpha(y) \right)$
		are local frame $a$-correspondents over $\algA$-frames.
		\item If $\phi$ and $\alpha$ are local frame $a$-correspondents,
		and $\psi$ and $\beta$ are local frame $a$-correspondents, then
		$\phi \wedge \psi$ and $\alpha \wedge \beta$ are local frame $a$-correspondents
		over $\algA$-frames.
		\item If $\phi$ and $\alpha$ are local frame $b$-correspondents for all $b \in A$,
		$\psi$ and $\beta$ are local frame $c$-correspondents for all $c \in A$ and
		$\phi$ and $\psi$ have no proposition variables in common, then
		$\phi \vee \psi$ and $\alpha \vee \beta$ are local frame $a$-correspondents
		over $\algA$-frames.
	\end{enumerate} 
\end{lemma}
\begin{proof}
    We prove the first and third claims and leave the second to the reader.  For    any $\algA$-frame $\modF$, 	
    \begin{align*}
		&\modF,w \Vdash_{a} \Box^{n} \phi\\
		\iff\quad &a \leq V(\Box^{n} \phi, w)=\bigwedge \{R^{n}wu \rightarrow V(\phi, u) \mid u\in W\}
		\text{ for all valuations }V\text{ on }\modF\\ 
		\iff\quad &a \leq R^{n}wu \rightarrow V(\phi, u)\text{ for all valuations }V\text{ on }\modF\text{ and }u\in W\\
		\iff\quad &a \wedge R^{n}wu  \leq V(\phi, u)\text{ for all valuations }V\text{ on }\modF\text{ and }u\in W\\
	\end{align*}
    where the first bi-implication is justified by Proposition \ref{Proposition:TruthPreserved} and Lemma \ref{lem:translate:boxed:atom}. Setting $b_u = a \wedge R^{n}wu$ for each $u \in W$, the last statement above is equivalent to  $\modF, u \Vdash_{b_u} \phi$ for each $u \in W$. This, by assumption, is equivalent to  $\modF \models_{b_u} \alpha(y) [y := u]$ for each $u \in W$. This, in turn, is equivalent to $a \wedge R^{n}wu \leq \Vert \alpha(y) [y := u] \Vert_{\modF}$, i.e., $a \leq R^{n}wu \to \Vert \alpha(y) [y := u] \Vert_{\modF}$, i.e., $a \leq  \Vert R^{n}xy \to \alpha(y) [x:= w, y := u] \Vert_{\modF}$ for all $u \in W$, which is equivalent to $a \leq  \Vert \forall y(R^{n}xy \to \alpha(y)) [x:= w] \Vert_{\modF}$, i.e., $\modF\models_{a} \forall y (R^{n}xy \rightarrow \alpha(y) [x:= w])$.
	
	Now consider (iii). By Lemma~\ref{Lemma:Disjunctions},
	\begin{align*}
		&\modF, w\Vdash_a\phi\vee\psi\\
		\iff\quad & \modF,w\Vdash_{a_1}\phi\text{ and }\modF,w\Vdash_{a_2}\psi
		\text{ for some }a_1,a_2\in A\text{ such that }a_1\vee a_2=a\\
		\iff\quad & \Vert\alpha [x := w]\Vert_{\modF}\geq a_1\text{ and }\Vert\beta [x:=w]\Vert_{\modF}\geq a_2
		\text{ for some }a_1,a_2\in A\text{ such that }a_1\vee a_2=a\\
		\iff\quad & \Vert\alpha\vee\beta\Vert_{\modF,v}=
		\Vert\alpha\Vert_{\modF,v}\vee\Vert\beta\Vert_{\modF,v}\geq a_1\vee a_2 =a\text{ for all }v\text{ on }\modF\\
		\iff\quad &\modF\models_a\alpha\vee\beta.
	\end{align*}
\end{proof}

\begin{thm}\label{thm:Sahl:Corresp:Classical:MV}
Let $a\in A$ and $\chi$ be a classical Sahlqvist formula. Let $\alpha$ be the local frame correspondent of $\chi$
obtained by the classical Sahlqvist-van Benthem algorithm.
Then $\chi$ and $\alpha$ are local frame $a$-correspondents over
$\algA$-frames.
\end{thm}

\begin{proof}
For definiteness' sake, we will take as our reference point the presentation of the Sahlqvist-van Benthem algorithm that can be found in \cite{BlueBook}. This algorithm finds correspondents for definite Sahlqvist implications. Correspondents for classical Sahlqvist formulas in general can than be found by applying (classical analogues) of Lemmas \ref{lem:Definite:Classical:Sahl} and \ref{Lemma:SahlFormulasToImplications}. Accordingly, it suffices for us to prove the
result for definite Sahlqvist implications.  Thus, suppose $\chi = \phi \rightarrow \psi$ is a definite classical Sahlqvist implication and let $\modF = (W,R)$ be any $\algA$-frame.
By Corollary~\ref{Corollary:StandardTranslation}, 
$\modF\Vdash_{a}\chi$ if,
and only if, 
$\modF \models_a \forall \oP\left(\ST_{x}(\chi)\right) [x:= w]$
where $\oP$ is a vector of all predicate variables corresponding to the propositional variables that occur in $\chi$.
(Note that since $\chi$ is in $\mvml$ and not $\emvml$,
the standard translation contains no occurrences of individual constants $\bc_\i$ and $\bc_\m$ or of nullary predicate variables $\bC_\i$ and $\bC_\m$.)

Recall that 
\begin{equation}\label{FOTransSahlqvist}
\forall \oP\left(\ST_{x}(\chi)\right)=
\forall \oP\left(ST_{x}(\phi) \rightarrow ST_{x}(\psi)\right)
\end{equation}
has $x$ as its only free variable and is clean, i.e., no quantifier binds $x$ and no two quantifiers bind the same variable.
Moreover, by the form of definite Sahlqvist antecedents, existential quantifiers in $\ST_{x}(\phi)$ are only in the scope of conjunctions and other existential quantifiers. Therefore, all existential quantifiers in the antecedent of \eqref{FOTransSahlqvist} can be moved to the 
front of the implication using the equivalences from 
Lemma~\ref{Lemma:Equivalences}, producing a formula of the form 
\begin{equation}\label{FOFormulaAfterEquivalences}
\forall \oP
\forall x_1\ldots \forall x_m(\Rel\wedge\Boxat\rightarrow \Pos),
\end{equation}
where  $\Pos$ is $\ST_x(\psi)$, $\Rel$ is a conjunction of first-order atomic formulas of the form 
$R x_ix_j$ coming from the translations of diamonds, and 
$\Boxat$ is a conjunction of translations of boxed atoms, i.e., formulas
of the form $\forall y\left(R^rx_iy \rightarrow P(y)\right)$.
Recall that if $r=0$, then $R^r x_i y$ is read as $x_i=y$.
Note that the free variables in $\Rel$ and $\Boxat$ are $x, x_1,\ldots, x_m$,
while only $x$ is free in $\Pos$.
Also observe that every occurrence of a unary predicate variable $P$ in the 
antecedent appears in $\Boxat$ in the
translation of a boxed atom, as described above.

Next we introduce two substitutions that produce instantiations:
the first is the instantiation which constitutes the minimal valuation in the classical case, while the second is minimal (in a sense that we clarify below) in the many-valued case. Let $\P$ be a unary predicate that occurs in (\ref{FOFormulaAfterEquivalences})
and let $\forall y_1\left(R^{r_1}x_{i_1}y_1 \rightarrow \P(y_1)\right),\ldots,
\forall y_k\left(R^{r_k}x_{i_k}y_k \rightarrow \P(y_k)\right)$
be all the translations of boxed atoms in $\Boxat$ in which $\P$ occurs.
Define a predicate substitution $\sigma$ such that, for each predicate variables $\P$, 
\[
\sigma\left(\P(y)\right):=
R^{r_{1}}x_{i_{1}}y\vee\cdots\vee R^{r_{k}}x_{i_{k}}y.
\]
Thus for a nullary predicate variable\footnote{Formally $\bC$ is $\bC_{\i}$ for some nominal $\i$ but, as the $\i$ plays no role, we will suppress the subscript.} $\bC$ we have
\[
\sigma_\bC\left(P(y)\right)=
\left(R^{r_{1}}x_{i_{1}}y\vee\cdots\vee R^{r_{k}}x_{i_{k}}y\right)\wedge \bC.
\]
In the classical case, $\sigma$ produces, set theoretically speaking, the smallest interpretation of the predicate symbols $\P$ which, relative to an assignment to $x_1\ldots x_m$ which makes $\Rel$ true, makes the antecedent $\Rel\wedge\Boxat$ true. Analogously, in the $\algA$-valued case, for any assignment $\overline{x}:=\overline{w}$ of $w, w_1, \ldots, w_m$ to $x, x_1,\ldots x_m$, and interpretation of the predicates $\obP:=\overline{P}$ that makes the antecedent non-zero (in $\algA$), there must by a $C \in \jty(\algA)$ such that $\Vert \sigma_\bC\left(\P(y)\right) [\overline{x}:=\overline{w}, y:= u, \bC := C]\Vert = \bigvee_{j=1}^k R^{r_j}w_{i_j}u\wedge C\leq P(u)$ for each predicate $\P$ and all $u\in W$ (for why this is so, please see the case (\ref{ThirdEquivalentStatement})$\Rightarrow$(\ref{FirstEquivalentStatement}), below). Thus, also here, $\sigma_\bC$ represents a `minimal' instantiation of the predicates. 

We write $\overline{x}$ for $x, x_1\ldots x_m$. Let $\overline{w}\in W$. To complete the proof we show that the following statements are equivalent:  
\begin{align}
&\modF\models_a\left(\Rel \wedge \Boxat \rightarrow \Pos\right)
[\overline{x}:=\overline{w},\obP:=\overline{P}] \quad \text{for all } \overline{P}:W\rightarrow \algA \label{FirstEquivalentStatement}\\
\iff\quad & \modF\models_a\left(\forall \bC\left(\Rel \wedge\sigma_\bC(\Boxat)
\rightarrow \sigma_\bC(\Pos)\right)\right)[\overline{x}:=\overline{w}]\label{SecondEquivalentStatement}\\
\iff\quad & \modF\models_a\left(\forall \bC\left(\Rel \wedge \bC \rightarrow 
\sigma_\bC(\Pos)\right)\right)[\overline{x}:=\overline{w}]\label{ThirdEquivalentStatement}\\
\iff\quad & \modF\models_a\left(\Rel \rightarrow \sigma(\Pos)\right)[\overline{x}:=\overline{w}]
\label{FourthEquivalentStatement}
\end{align}
\begin{enumerate}
\item[(\ref{FirstEquivalentStatement})$\Rightarrow$(\ref{SecondEquivalentStatement})] 
Since $\overline{P}$ is arbitrary, we can consider, in particular, the 
instantiation produced by the substitution $\sigma_\bC$.

\item[(\ref{SecondEquivalentStatement})$\Rightarrow$(\ref{ThirdEquivalentStatement})]
Suppose $\modF\models_a\left(\forall \bC\left(\Rel \wedge\sigma_\bC(\Boxat)
\rightarrow \sigma_\bC(\Pos)\right)\right)[\overline{x}:=\overline{w}]$.
Recall that $\Boxat$ is a conjunction of formulas of the form 
$\forall y\left(R^{r}x_{i}y \rightarrow \P(y)\right)$.
Then $\sigma_\bC(\Boxat)$ is a conjuction of formulas of the form 
\[
\forall y_\ell\left(R^{r_\ell}x_{i_\ell}y_\ell \rightarrow 
\left(\bigvee_{j=1}^k R^{r_{j}}x_{i_{j}}y_{\ell}\wedge \mathbf{C}\right)\right)
\]
where $R^{r_\ell}x_{i_\ell}y_\ell$ is the same as $R^{r_{j}}x_{i_{j}}y_{\ell}$ for some $1 \leq j \leq k$. This, together with Heyting algebra identity $a \to (a \vee b) = a \to b$, justifies the second equality below, while the inequality follows from the Heyting algebra inequality $c\leq(b\rightarrow c)$:
\begin{align*}
C &= \left\Vert\bC \right\Vert_{\modF}\\
&\leq \left\Vert\left(\forall y_\ell\left(R^{r_\ell}x_{i_\ell}y_\ell \rightarrow\bC\right)\right)
[\overline{x}:=\overline{w}]\right\Vert_{\modF}\\
&=\left\Vert\left(\forall y_\ell\left(R^{r_\ell}x_{i_\ell}y_\ell \rightarrow 
\left(\bigvee_{j=1}^k R^{r_{j}}x_{i_{j}}y_{\ell}\wedge \bC\right)\right)\right)
[\overline{x}:=\overline{w}]\right\Vert_{\modF}
\end{align*}
Thus, $C \leq \Vert\sigma_\bC(\Boxat)[\overline{x}:=\overline{w}]\Vert_{\modF}$.
%since each of its conjuncts is greater than or equal to $C$.
By the antitonicity of $\rightarrow$ in the first-coordinate, 
\[
\modF\models_a\left[\forall \bC\left(\Rel \wedge \bC \rightarrow 
\sigma_{\bC}(\Pos)\right)\right][\overline{x}:=\overline{w}].
\]
%%%%%% 13 ==> 11 %%%%%%%%%%%%%%%%%%%%
\item[(\ref{ThirdEquivalentStatement})$\Rightarrow$(\ref{FirstEquivalentStatement})]
We begin by proving the following chain of implications, where $C \in J^\infty(\algA)$:
\begin{align}
\begin{split}
C &\leq \Vert\left(\Rel\wedge\Boxat\right)[\overline{x}:=\overline{w},\obP:=\overline{P}]\Vert_\modF\\
&\Rightarrow \Vert\sigma_\bC(\P(y))[\overline{x}:=\overline{w},y:=u, \bC := C]\Vert_\modF \leq \Vert \P(y)[\obP:=\overline{P},y:=u]\Vert_\modF 
\quad \text{ for all }u\in W\\
&\Rightarrow 
\Vert\sigma_\bC(\Pos(x))
[\obP:=\overline{P},x:=u, \bC := C]\Vert_\modF
\leq \Vert\Pos[\obP:=\overline{P},x:=u]\Vert_\modF  \quad \text{ for all }u\in W
\end{split}\tag{$\dagger$}\label{eq:chain:for:13:i,plies:14}
\end{align}
Suppose $C \leq \Vert\left(\Rel\wedge\Boxat\right)[\overline{x}:=\overline{w},\obP:=\overline{P}]\Vert_\modF$
for some $C\in J^\infty(\algA)$, then $C \leq \Vert\Boxat[\overline{x}:=\overline{w}, \obP:=\overline{P}]\Vert_\modF$.\
Recall that $\Boxat$ is a conjunction of formulas of the form 
$\forall y\left(R^{r}x_{i}y \rightarrow \P(y)\right)$.
Hence, each conjunct satisfies 
$C \leq \Vert\forall y\left(R^{r}x_{i}y \rightarrow \P(y)\right)
[\overline{x}:=\overline{w},\obP:=\overline{P}]\Vert_\modF$.  Thus, for all $u\in W$ we have that
$C \leq R^{r}w_{i}u \rightarrow P(u)$, i.e.\ 
$R^{r}w_{i}u\wedge C\leq P(u)$.  But then 
$\bigvee_{j=1}^k R^{r_j}w_{i_j}u\wedge C\leq P(u)$ and it follows that, for all $u\in W$ and each $\P \in \obP$, 
\[
\Vert\sigma_\bC(\P(y))[\overline{x}:=\overline{w},y:=u, \bC := C]\Vert_\modF\leq
\Vert \P(y)[\obP:=\overline{P},y:=u]\Vert_\modF.
\]
Thus, my the monotonicity of positive formulas, this immediately implies that 
\[
\Vert\sigma_\bC(\Pos(x))[x:=w, \bC := C]\Vert_\modF\leq
\Vert\Pos[\obP:=\overline{P},x:=w]\Vert_\modF.
\]

Now let $S=\{C\in J^\infty(\algA)\mid C \leq \Vert\left(\Rel\wedge\Boxat\right)
[\overline{x}:=\overline{w},\obP:=\overline{P}]\Vert_\modF \}$.
Then,
\begin{eqnarray*}
&&\modF\models_a\left(\forall \bC\left(\Rel \wedge \bC \rightarrow 
\sigma_\bC(\Pos)\right)\right)[\overline{x}:=\overline{w}]\\
&\Rightarrow&\modF\models_a\left(\Rel \wedge \bC \rightarrow 
\sigma_\bC(\Pos)\right)[\overline{x}:=\overline{w}, \bC := C] \quad \text{for all } C \in S\\
&\Rightarrow&\modF\models_a \left(\Rel \wedge \bC \rightarrow 
\Pos\right)[\overline{x}:=\overline{w}, \obP:=\overline{P}, \bC := C] \quad \text{for all } C \in S.
\end{eqnarray*}
where the last implication follows by \eqref{eq:chain:for:13:i,plies:14}. 
Observe that the last line is equivalent to 
\[
a \leq \bigwedge_{C \in S} \Vert \left(\Rel \wedge \bC \rightarrow 
\Pos\right)[\overline{x}:=\overline{w}, \obP:= \overline{P}, \bC := C] \Vert_{\modF}
\]
and that $\bigvee S=\Vert\left(\Rel\wedge\Boxat\right)
[\overline{x}:=\overline{w},\obP:=\overline{P}]\Vert_\modF$. Using the complete distributivity of $\algA$ and the complete $\vee$-reversal of $\to$ in the first coordinate,  this implies that 
\[
a \leq \Vert \left(\Rel \wedge \Boxat \rightarrow 
\Pos\right)[\overline{x}:=\overline{w}, \obP:= \overline{P}] \Vert_{\modF},
\]
i.e., that
\[\modF\models_a\left(\Rel \wedge \Boxat \rightarrow \Pos\right)
[\overline{x}:=\overline{w},\obP:=\overline{P}].\]

\item[(\ref{ThirdEquivalentStatement})$\Rightarrow$(\ref{FourthEquivalentStatement})]
Suppose $\modF\models_a\left(\forall \bC\left(\Rel \wedge \bC \rightarrow 
\sigma_\bC(\Pos)\right)\right)[\overline{x}:=\overline{w}]$.
Since $\Pos$ is monotone, it follows from the defintions of $\sigma_\bC$ and $\sigma$
that $\left\Vert\sigma_\bC(\Pos)[\overline{x}:=\overline{w}]\right\Vert_{\modF}\leq
\left\Vert\sigma(\Pos)[\overline{x}:=\overline{w}]\right\Vert_{\modF}$.
But then $\modF\models_a\left(\forall \bC\left(\Rel \wedge \bC \rightarrow 
\sigma(\Pos)\right)\right)[\overline{x}:=\overline{w}]$ since $\rightarrow$
is order-preserving in the second coordinate. Then 
\begin{eqnarray*}
&&\modF\models_a\left(\forall \bC\left(\Rel \wedge \bC \rightarrow 
\sigma(\Pos)\right)\right)[\overline{x}:=\overline{w}]\\
&\iff& a \leq \bigwedge_{C\in J^{\infty}(\algA)}\Vert \left(\Rel \wedge \bC
\rightarrow \sigma(\Pos)\right)[\overline{x}:=\overline{w}, \bC := C] \Vert_{\modF}\\
&\iff& a \leq \Vert \left(\Rel \wedge \top 
\rightarrow \sigma(\Pos)\right)[\overline{x}:=\overline{w}] \Vert_{\modF}\\
&\iff& \modF\models_a\left(\Rel \rightarrow \sigma(\Pos)\right)[\overline{x}:=\overline{w}],
\end{eqnarray*}
where the second equivalence follows since  $\bigwedge_{c\in J^{\infty}(\algA)}(b \wedge c \to d) = (b \wedge \top) \to c$ in $\algA$.

\item[(\ref{FourthEquivalentStatement})$\Rightarrow$(\ref{ThirdEquivalentStatement})]
Suppose that $\modF\models_a\left(\Rel \rightarrow \sigma(\Pos)\right)[\overline{x}:=\overline{w}]$
and let $C\in J^\infty(\algA)$.  Then we have
\begin{align}
a &\leq \left\Vert\left(\Rel\rightarrow\sigma(\Pos)\right)[\overline{x}:=\overline{w}]\right\Vert_{\modF}\nonumber\\
&\leq \left\Vert\left(\Rel \wedge \bC \rightarrow\sigma(\Pos)\right)[\overline{x}:=\overline{w}]\right\Vert_{\modF}\label{eq:final:imp:2}\\
&= \left\Vert\left(\Rel \wedge \bC \rightarrow\sigma(\Pos)\wedge\bC\right)
[\overline{x}:=\overline{w}]\right\Vert_{\modF}\label{eq:final:imp:3}\\
&=
\left\Vert\left(\Rel \wedge \bC \rightarrow\sigma_\bC(\Pos)\wedge\bC\right)
[\overline{x}:=\overline{w}]\right\Vert_{\modF}\label{eq:final:imp:4}\\
&= \left\Vert\left(\Rel \wedge \bC \rightarrow\sigma_\bC(\Pos)\right)
[\overline{x}:=\overline{w}]\right\Vert_{\modF}\label{eq:final:imp:5},
\end{align}
where \eqref{eq:final:imp:2} and \eqref{eq:final:imp:5} follow by the antitonicity (monotonicity) of $\rightarrow$ in its first (second) coordinate, \eqref{eq:final:imp:3} is justified by the Heyting algebra identity $a \wedge c \to b = a \wedge c \to b \wedge c$, and \eqref{eq:final:imp:4} follows by Lemma~\ref{Lemma:EliminatingC}. Finally, as $C\in J^\infty(\algA)$ was arbitrary, we have that
\[
\modF\models_a\left(\forall \bC\left(\Rel \wedge \bC \rightarrow 
\sigma_\bC(\Pos)\right)\right)[\overline{x}:=\overline{w}].
\]
\end{enumerate}
\end{proof}

\section{Conclusions}\label{sec:conclusions}

We conclude by mentioning some areas for future research and further investigation.

\paragraph{Comparisons beyond Sahlqvist formulas.}  Example \ref{example:change:from:classical:corresp} suggests that the strong connection between classical and many-valued correspondents established in Theorem \ref{thm:Sahl:Corresp:Classical:MV} does not generalize to all inductive or even all Sahlqvist formulas, in the sense of Definition \ref{Sahlqvist:Form:Def}. Indeed, while Theorem \ref{thm:Sahl:Corresp:Classical:MV} establishes syntactic identity between classical and many-valued corespondents, Example \ref{example:change:from:classical:corresp} shows that, beyond the syntactic fragment of the classical Sahlqvist formulas, classical and many-valued corespondents need not even be equivalent when read as (crisp) first-order formulas. Any generalization of Theorem \ref{thm:Sahl:Corresp:Classical:MV} will have to take account of this, and it remains an open questions for which fragments of the inductive formulas such a generalizations are possible.

\paragraph{Other truth-value spaces.} Following \cite{Fitting, Fitting2} and a line of contributions (including \cite{Koutras}, \cite{Koutras1}, \cite{Koutras5}, \cite{Frank} and \cite{Frank1}) which develop correspondence-theoretic in this setting, we have focused on logics with Heyting algebras as truth-value spaces. The work in the present paper therefore also applies to the special case of G\"{o}del modal logic \cite{caicedo2010standard}. However, there is a large and growing body of work on many-valued modal logic based on more general truth-value algebras, e.g. the framework of \cite{Bou1} allows any resituated lattice. Several results in the present paper generalize to this setting in a reasonably straight forward way, including the standard translations, the definition of Sahlqvist and inductive formulas and the ALBA algorithm. However, results dealing with the syntactic shape of correspondents, including those of Section \ref{Section:SvBResult}, depend essentially on the fact that $\wedge$ and $\to$ form a residuated pair and therefore do not generalize in a transparent way. If and how this can be done remains an open question.

\bibliographystyle{siam}
\bibliography{M} 

\begin{thebibliography}{10}

\bibitem{Ackermann:Untersuchung}
{\sc W.~Ackermann}, {\em Untersuchung {\"{u}}ber das {E}liminationsproblem der mathematischen {L}ogic}, Mathematische Annalen, 110 (1935), pp.~390--413.

\bibitem{Hyb}
{\sc C.~Areces and B.~ten Cate}, {\em {H}ybrid {L}ogics}, in Handbook of Modal Logic, P.~Blackburn, J.~van Benthem, and F.~Elsevier, eds., vol.~3, Elsevier Science, 1~ed., 2007, ch.~3.

\bibitem{BadiaCaicedoNoguera2023}
{\sc G.~Badia, X.~Caicedo, and C.~Noguera}, {\em Frame definability in finitely valued modal logics}, Annals of Pure and Applied Logic, 174 (2023), p.~103273.

\bibitem{BlueBook}
{\sc P.~Blackburn, M.~de~Rijke, and Y.~Venema}, {\em Modal Logic}, Cambridge University Press, The Edinburgh Building, Cambridge, 2001.

\bibitem{Bou1}
{\sc F.~Bou, M.~Cerami, and F.~Esteva}, {\em {F}inite-valued {{\L}}ukasiewicz {M}odal {L}ogic is {P}{S}{P}{A}{C}{E}-complete}, in Proceedings of the Twenty-Second International Joint Conference on Artificial Intelligence - Volume Volume Two, AAAI Press, 2011, pp.~774--779.

\bibitem{caicedo2010standard}
{\sc X.~Caicedo and R.~O. Rodriguez}, {\em Standard g{\"o}del modal logics}, Studia Logica, 94 (2010), pp.~189--214.

\bibitem{chagrov1997modal}
{\sc A.~Chagrov and M.~Zakharyaschev}, {\em {M}odal {L}ogic}, Oxford logic guides, Clarendon Press, 1997.

\bibitem{conradie2024modal}
{\sc W.~Conradie, A.~De~Domenico, K.~Manoorkar, A.~Palmigiano, M.~Panettiere, D.~P. Prieto, and A.~Tzimoulis}, {\em Modal reduction principles across relational semantics}, Fuzzy Sets and Systems, 481 (2024), p.~108892.

\bibitem{CoGhPa14}
{\sc W.~Conradie, S.~Ghilardi, and A.~Palmigiano}, {\em Unified correspondence}, in Johan van Benthem on Logic and Information Dynamics, A.~Baltag and S.~Smets, eds., vol.~5 of Outstanding Contributions to Logic, Springer International Publishing, 2014, pp.~933--975.

\bibitem{ALBA}
{\sc W.~Conradie and A.~Palmigiano}, {\em {A}lgorithmic {C}orrespondence and {C}anonicity for {D}istributive {M}odal {L}ogic}, Annals of Pure and Applied Logic, 163 (2012), pp.~338 -- 376.

\bibitem{nonDistALBA}
{\sc W.~Conradie and A.~Palmigiano}, {\em Algorithmic correspondence and canonicity for non-distributive logics}, Annals of Pure and Applied Logic, 170 (2019), pp.~923--974.

\bibitem{conradie2017algebraic}
{\sc W.~Conradie, A.~Palmigiano, and S.~Sourabh}, {\em Algebraic modal correspondence: Sahlqvist and beyond}, Journal of Logical and Algebraic Methods in Programming, 91 (2017), pp.~60--84.

\bibitem{DGP}
{\sc J.~M. Dunn, M.~Gehrke, and A.~Palmigiano}, {\em Canonical extensions and relational completeness of some substructural logics}, {J}ournal of {S}ymbolic {L}ogic, 70(3) (2005), pp.~713 -- 740.

\bibitem{Koutras}
{\sc P.~E. Eleftheriou and C.~D. Koutras}, {\em {F}rame {C}onstructions, {T}ruth {I}nvariance and {V}alidity {P}reservation in {M}any-valued {M}odal {L}ogic}, Journal of Applied Non-Classical Logics, 15 (2005), pp.~367 -- 388.

\bibitem{Fitting}
{\sc M.~C. Fitting}, {\em {M}any-valued {M}odal {L}ogic}, Fundamenta Informaticae, 15 (1991), pp.~235 -- 254.

\bibitem{Fitting2}
\leavevmode\vrule height 2pt depth -1.6pt width 23pt, {\em Many-valued modal logics 2}, Fundamenta Informaticae, 17 (1992), pp.~55 -- 73.

\bibitem{Frank}
{\sc S.~Frankowski}, {\em {D}efinable {C}lasses of {M}any-{V}alued {K}ripke {F}rames}, Bulletin of the Section of Logic, 35 (2006), pp.~27 -- 36.

\bibitem{Frank1}
\leavevmode\vrule height 2pt depth -1.6pt width 23pt, {\em {G}eneral {A}pproach to {M}any-valued {K}ripke {M}odels}, Bulletin of the Section of Logic, 35 (2006), pp.~11--26.

\bibitem{GEHRKE200565}
{\sc M.~Gehrke, H.~Nagahashi, and Y.~Venema}, {\em {A} {S}ahlqvist {T}heorem for {D}istributive {M}odal {L}ogic}, Annals of Pure and Applied Logic, 131 (2005), pp.~65 -- 102.

\bibitem{Goranko2006-GORECF}
{\sc V.~Goranko and D.~Vakarelov}, {\em {E}lementary {C}anonical {F}ormulae: {E}xtending {S}ahlqvist's {T}heorem}, Annals of Pure and Applied Logic, 141 (2006), pp.~180--217.

\bibitem{hajek1998metamathematics}
{\sc P.~H{\'a}jek}, {\em Metamathematics of Fuzzy Logic}, Trends in Logic, Springer Netherlands, 1998.

\bibitem{Koutras5}
{\sc C.~D. Koutras}, {\em {A} {C}atalog of {W}eak {M}any-valued {M}odal {A}xioms and their {C}orresponding {F}rame {C}lasses}, Journal of Applied Non-Classical Logics, 13 (2003), pp.~47 -- 71.

\bibitem{kracht1999tools}
{\sc M.~Kracht}, {\em Tools and techniques in modal logic}, vol.~97, Elsevier Amsterdam, 1999.

\bibitem{Koutras1}
{\sc C.~Nomikos, C.~D. Koutras, and P.~Peppas}, {\em {C}anonicity and {C}ompleteness {R}esults for {M}any-valued {M}odal {L}ogics}, Journal of Applied Non-Classical Logics, 12 (2002), pp.~7 -- 41.

\bibitem{Sahl}
{\sc H.~Sahlqvist}, {\em {C}orrespondence and {C}ompleteness in the {F}irst- and {S}econd-order {S}emantics for {M}odal {L}ogic}, in Proceedings of the Third Scandinavian Logic Symposium, 1975.

\bibitem{schotch1978note}
{\sc P.~K. Schotch, J.~B. Jensen, P.~F. Larsen, and E.~J. MacLellan}, {\em A note on three-valued modal logic}, Notre Dame Journal of Formal Logic, 19 (1978), pp.~63--68.

\bibitem{SEGERBERG}
{\sc K.~Segerberg}, {\em {S}ome {M}odal {L}ogics {B}ased on a {T}hree-valued {L}ogic}, Theoria, 33 (1967), pp.~53--71.

\bibitem{thomason1978possible}
{\sc S.~Thomason}, {\em Possible worlds and many truth values}, Studia Logica: An International Journal for Symbolic Logic, 37 (1978), pp.~195--204.

\bibitem{vanb}
{\sc J.~van Benthem}, {\em Modal Logic and Classical Logic}, Bibliopolis, Napoli, 1983.

\bibitem{vanBenthemChapter2001}
\leavevmode\vrule height 2pt depth -1.6pt width 23pt, {\em {C}orrespondence {T}heory}, in Handbook of Philosophical Logic, D.~Gabbay and F.~Guenthner, eds., Springer, Heidelberg, 2001, ch.~4, pp.~325--408.

\bibitem{van2005minimal}
{\sc J.~Van~Benthem}, {\em Minimal predicates, fixed-points, and definability}, The Journal of Symbolic Logic, 70 (2005), pp.~696--712.

\end{thebibliography}

\end{document}